\documentclass[pdflatex, sn-mathphys-num]{sn-jnl}

\usepackage{xcolor}

\usepackage[T1]{fontenc}
\usepackage{amsmath}
\usepackage[title]{appendix}
\usepackage{amsfonts}
\usepackage{bm}
\usepackage{a4wide}
\usepackage{hyperref}
\usepackage{subcaption}
\usepackage{graphicx}
\usepackage{enumitem}
\usepackage[ruled,vlined]{algorithm2e}
\SetInd{5pt}{10pt}
\SetNlSty{relax}{}{\enspace}
\usepackage[standard]{ntheorem}
\definecolor{mycustomcolor}{RGB}{128, 0, 128}

\DeclareMathOperator*{\argmin}{arg\,min}
\newcommand{\norm}[1]{\left\lVert #1 \right\rVert}
\newcommand{\nens}{N}

\newcommand{\statebis}{\bm{z}}
\newcommand{\bx}{\bm{x}}

\newcommand{\vU}{\bm{\mathrm{U}}}
\newcommand{\vu}{\bm{\mathrm{u}}}

\newcommand{\Fcorr}{\bm{F}}

\newcommand{\sigmaZm}{0.5}
\newcommand{\meanZm}{\pi/2 + 0.6}
\newcommand{\smLow}{0.8}
\newcommand{\smUp}{1.2}

\newcommand{\visc}{D}

\newcommand{\cN}{\mathcal N}
\newcommand{\cU}{\mathcal U}

\newcommand{\refv}{1.0}
\newcommand{\refvisc}{0.05}
\newcommand{\Dlow}{0.02}
\newcommand{\Dup}{0.08}
\newcommand{\vmean}{0.9}
\newcommand{\vstd}{1.2}

\newcommand{\Cov}{\bm{P}}
\newcommand{\sigz}{\sigma_0^2 = 0.5}
\newcommand{\zz}{x_0 = 0.02}

\begin{document}

\title[Ensemble Data Assimilation for Particle-based Methods]{Ensemble Data Assimilation for Particle-based Methods}
\author*[1]{\fnm{Marius} \sur{Duvillard}}\email{marius.duvillard@gmail.com}
\author[1]{\fnm{Loïc} \sur{Giraldi}}\email{loicgiraldi@gmail.com}
\author[2]{\fnm{Olivier} \sur{Le Maître}}\email{olivier.le-maitre@polytechnique.edu}
\affil*[1]{\orgdiv{DES, IRESNE, DEC, SESC, LMCP}, \orgname{CEA}, \orgaddress{\street{Cadarache},  \postcode{F-13108}, \city{Saint-Paul-Lez-Durance}, \country{France}}}
\affil[2]{\orgdiv{Centre de Mathématiques Appliquées}, \orgname{Ecole Polytechnique, IPP}, \orgaddress{\street{Route de Saclay}, \postcode{91128}, \city{Palaiseau}, \country{France}}}
\date{}


\abstract{
    This study presents a novel approach to applying data assimilation techniques for particle-based simulations using the Ensemble Kalman Filter. While data assimilation methods have been effectively applied to Eulerian simulations, their application in Lagrangian solution discretizations has not been properly explored. We introduce two specific methodologies to address this gap. The first methodology employs an intermediary Eulerian transformation that combines a projection with a remeshing process. The second is a purely Lagrangian scheme that is applicable when remeshing is not adapted. These methods are evaluated using a one-dimensional advection-diffusion model with periodic boundaries. Performance benchmarks for the one-dimensional scenario are conducted against a grid-based assimilation filter. Subsequently, assimilation schemes are applied to a non-linear two-dimensional incompressible flow problem, solved via the Vortex-In-Cell method. The results demonstrate the feasibility of applying these methods in more complex scenarios, highlighting their effectiveness in both the one-dimensional and two-dimensional contexts.}


\maketitle


\section{Introduction}
Numerical simulation enables the assessment of complex real-world systems, for instance, to facilitate the optimization of complex systems and perform risk analysis at reduced costs.
With the advancement of computational resources, simulation has become an essential tool for understanding and designing processes, particularly in mechanical engineering.
Simulation of mechanical systems has historically used grid-based or mesh-based methods. These techniques, which use structured or unstructured meshes, are known for their convergence properties and computational efficiency. Meshless methods offer significant promises for complex physics with large deformations (moving interfaces, fragmentation, or distortion).

Meshless methods, specifically particle-based methods, discretize the computational domain by a collection of particles (or computational elements) that move with the dynamic in a Lagrangian fashion. Each particle transports material properties and internal variables. Particles can discretize a continuum medium and are associated with a kernel to reconstruct continuous fields and differential operators. In this work, we will mainly focus on the Vortex Method~\cite{cottet_vortex_2000,mimeau_review_2021}, which discretizes the vorticity field of the flow of an incompressible fluid using the stream function-vorticity formulation of the Navier-Stokes equations.

When simulating a physical system, the numerical solution contains errors that need to be understood, quantified, and minimized. When observations are available, integrating this information can lead to a more accurate estimation of the system state. In this context, data assimilation techniques aim at combining various sources of information to improve the estimation of the state. Integrating model predictions and observational data has been widely applied in disciplines such as meteorology, oceanography, hydrology, and geosciences~\cite{bocquet_introduction_2014}.

In data assimilation, two prominent families of approaches have emerged: variational and stochastic methods. Variational approaches~\cite{variational_method} focus on finding an optimal system state by minimizing a cost function that measures the misfit between model predictions and observations while penalizing the magnitude for the correction from the prior state. The most common formulations derive from 3DVar and 4DVar~\cite{talagrand1997assimilation}.

The stochastic approaches go beyond mere state estimation; they delve into the quantification of uncertainty associated with the estimated states. Uncertainty quantification is a critical aspect, especially in dynamic systems, for reliable decision-making and model improvement. The assimilation process is performed within a Bayesian framework with a forecast and an analysis step. The Kalman filter~\cite{kalman_new_1960} is a sequential Bayesian filter based on a linear model and Gaussian distribution assumptions. More advanced filters handle nonlinear models and arbitrary state distributions. One of the most popular Bayesian filters is the Ensemble Kalman Filter (EnKF), introduced in~\cite{evensen_sequential_1994}. Its popularity derives primarily from its capacity to deal with high-dimensional problems and its remarkable robustness when deviating from the initial Gaussian assumptions. It consists of approximating the first two moments of the probability distribution of a state thanks to an ensemble of simulations, called members.

EnKF has been extensively employed in Eulerian discretization frameworks, i.e., simulation over a fixed computational mesh. However, its application to Lagrangian approaches presents unique challenges, especially when the Lagrangian discretizations are different and adapted for each member of the ensemble. Each member could involve different particle sets differing in position and size. A first solution consists of defining the state considering the position and associated intensities of each particle~\cite{darakananda_data-assimilated_2018,chen_superfloe_2022}. However, it suggests a similar number of particles with unambiguous indexation between members; in the other case, it leads to an exponential increase of particles by using an extend-and-pad approach. 

Eulerian methods face similar difficulties in the context of Multi-Resolution Analysis (MRA) or moving mesh simulations. It has been proposed to introduce a suitable approximation space (with projection, interpolation, and restriction operators) to represent the members and the corrections~\cite{siripatana_combining_2019,bonan_data_2017}. Different choices of the approximation spaces yield a spectrum of methods.

In the Lagrangian case, a solution consists of considering an arbitrary reference discretization space for all members, for instance, using a fixed grid of particles. This way, the member's state can be defined on the same approximation space using projection and interpolation operators.

However, for specific Lagrangian methods, the use of a regridding operator is infeasible. More specifically, the change of the particle positions is subjected to constraints. Consequently, our objective is to address this limitation by developing a formulation that preserves the particle-based representation of each member by correcting the intensities only.

This paper introduces new approaches to apply ensemble data assimilation techniques to meshless simulations that discretize a continuum domain. First, the Remesh-EnKF uses a new reference particle discretization. This way, the state is updated with a controlled number of particles. This first method is based on the complete regridding of the members' particle discretization on which the classical EnKF analysis is performed.
Then, if the particle discretization needs to be preserved, the Part-EnKF is introduced. In this case, the analyzed field is approximated by each member's previous particle discretization. The particle positions are unchanged; only the strengths are modified by regression or approximation of the analyzed solution.

In the following, background on data assimilation and the EnKF correction are introduced in Section~\ref{Background_DA}, as well as particle-based methods in Section~\ref{Background_Part}. Then, alternative methods are described in Section~\ref{Methods}. Afterward, those filters are compared with a grid-based filter in a 1D Advection-Diffusion problem in Section~\ref{App_1D}, and an incompressible viscous flow is solved using a Vortex Method in Section~\ref{App_2D} where the filters are quantitatively analyzed. Major conclusions from this work are provided in Section~\ref{sec:conclusion_chap4}.
\section{Background methods}
\subsection{Data assimilation}\label{Background_DA}
Data assimilation is the use of observations to correct a system's predicted state. It is generally formulated within a probabilistic framework. It allows for rigorously dealing with observations and model errors to improve the estimate of the real state and associated uncertainty. To this end, state and observations are modeled as random variables. It is achieved by establishing a recurrence in probability distributions through Bayesian estimation to determine the posterior of the state given the current observation.

\subsubsection{Data assimilation setting}
A hidden Markov chain is used to model the evolution of state and observation. We position ourselves within a finite-dimensional state. The forecast and observation are introduced, such that for $ k \geq 0$,
\begin{equation*}
    \begin{cases}
        \statebis_{k+1} = \mathcal{M}_{k+1} (\statebis_{k}) + \bm{\eta}_{k+1}, \\
        \bm{y}_{k+1} = \mathcal{H}_{k+1} (\statebis_{k+1}) + \bm \varepsilon_{k+1},
    \end{cases}
\end{equation*}where $\mathcal{M}_{k+1}$ is the model operator that maps the state between time $k$ and $k+1$, and $\mathcal{H}_k$ is the observation operator. The vector $\statebis_k \in \mathbb{R}^n$ is the state at time $k$, $\bm{y}_k \in \mathbb{R}^m$ the observation vector, $\bm{\eta}_{k}$ is the model error that accounts for discrepancies between the physics and its model, and $\bm{\varepsilon}_k$ is the observation error which combines measurement error and representation error. We assume that $\bm{\eta}_{k}, \bm{\varepsilon}_k$ are random variables following Gaussian distributions with zero mean and covariance matrices $\bm Q_k$ and $\bm R_k$ respectively. Finally,$\statebis_k \perp \bm \varepsilon_{k' \geq k}$ as well as $\statebis_k \perp \bm \eta_{k' \geq k}$. 

\subsubsection{Ensemble Kalman Filter}\label{enkf}
The filtering problem consists of assessing the current state distribution by utilizing observed data up to the present time. The Markovian structure allows the assimilation to be performed sequentially following a forecast step and an analysis step after the acquisition of a new observation. In the rest of the section, we remove the time index $k$ and present the forecast and analysis steps for a single assimilation step.

The Kalman filter~\cite{kalman_new_1960} is a Bayesian filter that, in addition to the previously mentioned assumptions, requires $\mathcal{M}$ and $\mathcal{H}$ to be linear operators and the initial state Gaussian. In this case, the posterior distribution of the state is also Gaussian and fully characterized by its mean and covariance. The Kalman estimator is thus a recursive version of the Minimum Mean Square Error applied to the Gaussian linear model.

The analysis is based on the combination of the forecast state, defined by the mean and covariance $(\bm z^f, \bm P)$, with observation $\bm y$ based on the Kalman gain $\bm K$ such that the analysis state $\bm z^a$ is
\begin{equation}\label{eq:kalman}
    \bm z^a = \bm z^f + \bm K (\bm y - \bm H z^f), \quad \text{with } \bm K = \bm P \bm H^T(\bm H \bm P \bm H^T + \bm R)^{-1}.
\end{equation}

The Ensemble Kalman Filter (EnKF) introduced in~\cite{evensen_sequential_1994} is a data assimilation method adapted to high dimensional non-linear problems. The distribution of the state is still described by its first two moments, which are estimated from a finite ensemble of states.

We present the stochastic EnKF, where the observations are perturbed to account for observation errors. Assuming we have an ensemble of $N$ members $\left\{\bm \statebis_i \right\}_{i=1}^N$, we forecast the ensemble by propagating each state with the dynamic model and obtain a forecast \mbox{ensemble: $\{\bm z_i^f = \mathcal M(\bm z_i)\}_{i=1}^N$}.

The Kalman gain $\bm K$ is defined by approximating covariance matrices using empirical estimations:
\begin{eqnarray*}
    \Cov \bm H^T &\simeq& \frac{1}{N - 1} \sum_{i = 1}^{N} {(\statebis_i^f - \overline{\statebis}^f)} {\left[ \mathcal{H}(\statebis_i^f) - \overline{\bm{y}}^f\right]}^T \\
    \bm H \Cov \bm H^T &\simeq& \frac{1}{N -1} \sum_{i = 1}^{N}\left[ \mathcal{H}(\statebis_i^f) - \overline{\bm{y}}^f\right] {\left[ \mathcal{H}(\statebis_i^f) - \overline{\bm{y}}^f\right]}^T \\
    \text{where } \overline{\statebis}^f &=& \frac1N \sum_{i=1}^N \statebis_i^f, \quad \text{and } \overline{\bm{y}}^f = \frac1N \sum_{i=1}^N \mathcal{H}(\statebis_i^f)
\end{eqnarray*}

The Kalman update in Equation~\eqref{eq:kalman} is then applied to each member $\bm z^f_i$, leading to a new analysis ensemble $\{\bm z^a_i\}_{i=1}^N$.
The EnKF updates can be defined in the \mbox{members'} space such that
\begin{eqnarray}~\label{enkf_formula_Fcorr}
    \statebis^a_i & = & \statebis_i^f + \frac{1}{N-1}\sum_{j=1}^N \left[(\mathcal H(\statebis_j^f) - \overline{\bm{y}})^T (\bm C_H + \bm R)^{-1} (\bm y + \bm \varepsilon_{i} - \mathcal H(\statebis_i^f))\right] \statebis_j^f \nonumber \\
    & = & \statebis_i^f + \sum_{j=1}^N F_{ji} \statebis_j^f,
\end{eqnarray}where the matrix of coefficients $\Fcorr$ only depends on the ensemble members through the predicted observations ensemble $\left\{\mathcal{H}(\statebis_i^f) \right\}^N_{i=1}$, the observation $\bm y$ and the associate perturbations~$\bm \varepsilon_i \sim \mathcal N(\bm 0, \bm R)$.
\subsection{Particle methods}\label{Background_Part}
Particle methods are used to solve continuous problems in fluid and solid mechanics. We detail our Lagrangian assimilation method for the particle discretization of a field $u$ on a domain $\Omega \subseteq \mathbb{R}^d$.

A real valued function $u$ defined on $\mathbb R^d$ can be written as
\begin{equation*}
    u(\bm x) = \int u(\bm x') \delta(\bm x' - \bm x) d\bm x',
\end{equation*}with $\delta$ the Dirac delta function. Introducing the kernel function $\phi_\varepsilon$ defined by
\begin{equation*}
    \phi_\varepsilon(\bm x) = \frac{1}{\varepsilon^d} \phi\left(\bm x / \varepsilon\right),
\end{equation*}and such that
\begin{equation*}
    \int \phi(\bm x) d\bm x = 1, \quad \lim_{\varepsilon \to 0} \phi_\varepsilon(\bm x) = \delta(\bm x).
\end{equation*}We get the smooth version $u_\varepsilon$ of $u$ as
\begin{equation*}
    u_\varepsilon(\bm x) = \int u(\bm x') \phi_\varepsilon(\bm x-\bm x') d\bm x.
\end{equation*}
In the kernel definition, $\varepsilon > 0$ is called the smoothing length. We discretize the convolution with a set of $N_p$ particles $\mathcal P \doteq \left\{\bm x_p, U_p\right\}_{p = 1}^{N_p}$, where $\bm x_p \in \mathbb R^d$ are the spatial coordinates of the $p$-th particle and $U_p \in \mathbb{R}$ its intensity. The integral is approximated as
\begin{equation}~\label{part_approx}
    u(\bm x) \simeq \sum_{p \in \mathcal P} U_p \phi_\varepsilon (\bm x-\bm x_p).
\end{equation}

We present several ways to construct the particle approximation of the continuous field, i.e., determine the strengths $U_p$, given a set of position $\bm x_p$.
The particle approximation can readily be extended to $u(x) \in \mathbb{R}^{n > 1}$.

\subsubsection{Construction of particle approximation}\label{interpOp}
A first particle approximation uses particle intensities defined by
\begin{equation}\label{eq:first_part_approx}
    U_p = \int_{\Omega_p} u(\bm x) d\bm x \simeq u(\bm x_p)~V_p,
\end{equation}where $V_p$ is the volume of the $p$-th particle. A better approximation is obtained using Beale's formula~\cite{beale_accuracy_1988}, which correct iteratively the intensities $U_p$ in order to recover the value $u(\bm x)$ at the particle locations. For $k \geq 0$, it gives iterative values until convergence:
\begin{equation*}
    U^{k+1}_p = U^{k}_p + \left(u(\bm x_p) - \sum_{p' \in \mathcal P} U^k_{p'} \phi_{\varepsilon}(\bm x_{p'} - \bm x_{p})\right)~V_p.
\end{equation*}

Another approach consists of using regression methods. For a given set of particle positions, the vector of intensities $\bm{\mathrm{U}} = [U_1, \dots, U_{N_p}]^T$ is obtained by minimizing the approximation error
\begin{equation*}
    e^2 = \sum_{p' \in \mathcal P} | u(\bm x_{p'}) - \sum_p U_p \phi_\varepsilon (\bm x_{p'} - \bm x_p)|^2.
\end{equation*}

We define $\bm{\mathrm{u}} = [ u(\bm x_1), \dots,u(\bm x_{N_p})]^T$ the vector of values of $u$ at the particle location and $\bm \Phi$ the matrix of kernel where $\Phi_{ij} = \phi_\varepsilon(\bx_i - \bx_j)$. The minimisation of $e^2$ leads to the following problem:
\begin{equation}\label{eq:min_prob_intensities}
    \bm{\mathrm{U}}^*= \argmin_{\bm{\mathrm{U}} \in \mathbb R^{N_p}} \norm{\bm{\mathrm{u}} - \bm \Phi \bm{\mathrm{U}}}^2_2,
\end{equation}where the solution is $\bm{\mathrm{U}}^* = (\bm \Phi^T \bm \Phi)^{-1} \bm \Phi \bm{\mathrm{u}}$. The problem can be ill-posed, particularly for non-well-distributed particles. Problem~\eqref{eq:min_prob_intensities} can be regularized by introducing a penalization term. The Ridge regression introduces a penalization of the form $\lambda \norm{\bm{\mathrm{U}}}_2^2$, where $\lambda \geq 0$, to get
\begin{equation*}
    \bm{\mathrm{U}}_{\text{ridge}}^* = \argmin_{\bm{\mathrm{U}} \in \mathbb R^{N_p}} \norm{\bm{\mathrm{u}} - \bm \Phi \bm{\mathrm{U}}}_2^2 + \lambda \norm{\bm{\mathrm{U}}}^2_2,
\end{equation*}with the following solution $\bm{\mathrm{U}}^*_{\text{ridge}} = (\bm \Phi^T \bm \Phi + \lambda \bm{\mathrm{I}})^{-1} \bm \Phi \bm{\mathrm{u}}$.

These approach requires the resolution of a $N_p$-by-$N_p$ system and becomes infeasible for large particle sets. Beale's formula is therefore preferred, although it lacks regularization and may present convergence issues.

\subsubsection{Remeshing particle discretization}\label{remesh_part}
To minimize distortions of the particle distribution and improve the particle approximation, a new particle discretization is often generated~\cite{cottet_vortex_2000,cottet_multi-purpose_1999}.
Typically, the remeshing creates a new set of particles associated with a uniform grid with step size $d_p$, the characteristic size of the particles.
The process is designed to preserve the first moments of the original approximation. In our work, we propose a two-step approach. First, we execute a projection step~\ref{it:assigment_art} to transfer the particle intensities to a grid. Subsequently, an interpolation step~\ref{it:interpolation_art} is performed to yield a new set of regularly spaced particles.

Our presentation is restricted to the one-dimensional case, where $\Omega \subseteq \mathbb R$. The extension in higher dimensions is achieved through the tensorization of the one-dimensional formulas and operators.
\begin{enumerate}[label=(\alph*)]
    \item \textbf{Projection on an Eulerian grid}\label{it:assigment_art}

          The nodes are defined on a regular grid with spacing $\ell = 2 d_p$. We denote by $x^g_{i} = i~\ell, \; i \in \mathbb Z$ the grids point location, and the node values as $\mathrm u^g_i$ with an associate volume $V_i$ (typically  $V_i = \ell^d$). By using suitable shape function $W$, the projection step distributes the particles intensity to the nodes:
          \begin{equation}\label{eq:nodal_values}
              \mathrm u^g_i = \frac1{V_i} \sum_{p \in \mathcal P} U_p W \left(\frac{x_i^g - x_p}{\ell} \right), \; \forall i \in \mathbb Z.
          \end{equation}The function $W$ determines a redistribution of the intensity on the grid and should at least verified:
          \begin{equation*}
              \sum_{i \in \mathbb Z} W \left(\frac{x_i^g - x}{\ell} \right) = 1, \quad \forall x \in \Omega.
          \end{equation*}
    \item \textbf{Interpolation on a new regular particle discretization}\label{it:interpolation_art}
          The nodal values yield the grid approximation $u^g$:
          \begin{equation}\label{eq:grid_approx}
              u^g(x) = \sum_{i \in \mathbb Z}\mathrm u^g_i W \left(\frac{x - x^g_i}{\ell} \right) \quad \forall x \in \Omega.
          \end{equation}
          A new set of particles is defined at $1/4$ and $3/4$ of each cell between nodes: $x^*_{p} = p~d_p - d_p/2, \; p \in \mathbb Z$. For this new particle set, we define the intensities $U^*_{p}$ using the first approximation in Equation~\eqref{eq:first_part_approx}
          \begin{equation}\label{eq:new_set_intensities}
              U^*_{p} = u^g(x^*_{p}) V_{p} = V_{p} \sum_{i \in \mathbb Z}\mathrm u^g_i W \left(\frac{x^*_{p} - x^g_i}{\ell} \right),
          \end{equation}and the field $u^g$ defined in Equation~\eqref{eq:grid_approx} leading to the following solution
          \begin{equation}\label{eq:interpolation_part}
              u^*(\bx) = \sum_{p\in\mathcal P^*} U^*_{p} \phi_\varepsilon(x - x^*_{p}).
          \end{equation}
\end{enumerate}

The combination of both steps generates a new undistorted particle distribution. The function $W$ determines the number of moments conserved. 
The piecewise linear interpolation function ensures the conservation of moment 0. More involved construction of function $W$ is needed for the conservation of higher moments. While B-splines improve the smoothness, they are only exact for linear functions. \cite{monaghan_extrapolating_1985} proposes a systematic approach to enhance accuracy. The concept involves constructing a new shape function based on a cutoff and its radial derivative. For instance, the cubic B-spline is improve up to be exact for quadratic functions, by using the following interpolation kernel
\begin{eqnarray*}~\label{cubic_radial_kernel}
    M_4'(r) &=& \left\{ \begin{aligned}
         & 1 - \frac{5}{2}r^2 + \frac{3}{2} |r|^3 & 0 \leq & |r| \leq 1 & \\
         & \frac{1}{2}{(2 - |r|)}^2(1 - |r|)      & 1 \leq & |r| \leq 2 & \\
         & 0                                      & 2 \leq & |r|.
    \end{aligned}
    \right.
\end{eqnarray*}
One interest in using such kernels is that they involve compact stencils, i.e. particles are distributed over a limited set of nodes. Ensuring higher moments conservation comes at the cost of losing positivity of $W$. In this work, we utilize the $M_4'$ kernel. Finally, in dimension $d > 1$, the redistribution kernels $W$ are obtained as the product of the one-dimensional kernel along each dimension
\begin{equation*}
    W_d(\bm r) = \prod_{i = 1}^d W (r_{i}).
\end{equation*}
\section{Assimilation for Lagrangian simulations}\label{Methods}
\subsection{Challenges}
This section outlines the development of ensemble data assimilation techniques tailored for particle-based simulations. While the forward step of EnKF is classically performed, the analysis is defined in Equation~\eqref{enkf_formula_Fcorr}. It consists of a correction of the members' state by linear combinations of the forward states. The coefficient of the correction, given by the matrix $\Fcorr$, involves only the observations, predictions, and noise. The matrix $\Fcorr$ can thus be evaluated without particular issues.

In contrast, forming the correction by linearly combining the forward states is not trivial in the Lagrangian context, especially when the members use different particle sets. This situation prevents the direct linear combination of the vector of intensities, whose components are associated with different locations and may have different sizes.

One possibility, as explored in~\cite{darakananda_data-assimilated_2018,le_provost_ensemble_2021}, is to extend consistently the particle sets of each member to the union of all the members particle set, $\mathcal{P}^a = \bigcup_k \mathcal{P}_k^f$. In this increased discretization space, the states write as
\begin{equation}\label{eq:intensities_expanded}
    u_i^f(\bx) = \sum_{p \in \mathcal P^a} \tilde{U}_{i,p} \phi_\varepsilon(\bx - \bx_p), \quad \text{where } \tilde U^f_{i,p} = \left\{\begin{aligned}
         & U^f_{i,p} \; & \text{if } p \in \mathcal P_i^f \\
         & 0   \;       & \text{otherwise}
    \end{aligned}\right.
\end{equation}

By setting $\tilde{\vU}^f_{i} \in \mathbb R^{|\mathcal P^a|}$ the vector of intensities, they are updated to
\begin{equation}
    \vU_i^a = \tilde{\vU}^f_i + \sum_{j = 1}^N F_{ji} \tilde{\vU}_j^f.
\end{equation}

In general, the zero intensities of $\tilde{\vU}^f_{i}$ in~\eqref{eq:intensities_expanded} will not remain zero in $\vU_i^a$, preventing its reduction. Consequently, this approach induces an exponential growth of the particle discretization. Therefore, one has to use this extend-and-pad approach with a reduction strategy that removes particles from $\mathcal P^a$ depending on the member. This can involve thresholding, aggregation, and remeshing~\cite{darakananda_data-assimilated_2018, yue_continuum_2015,cottet_multi-purpose_1999}.

To avoid having to deal with large particle sets and designing a subsequent reduction strategy, we apply a more direct treatment. Starting from the continuous approximations of the states, we write
\begin{equation}~\label{eq:analysed_field_F}
    u^a_i(\bx) = u^f_i(\bx) + \sum_{j=1}^N F_{ji} u^f_j(\bx) \quad i = 1,\dots, \nens.
\end{equation}

At this point, the question amounts to the construction of particle approximations of the $u_i^a$. We propose two approaches.
\begin{itemize}
    \item The Remesh-EnKF filter, in Section~\ref{sec:remesh_enkf}, employs a remeshing technique to come up with a unique particle sets $\mathcal P^g$ for discretizing all members. Consequently, classical EnKF can be applied to the intensities through standard linear combinations while controlling the total number of particles.
    \item The Part-EnKF filter, in Section~\ref{sec:method_part_enkf}, aims at reconstructing an approximation of $u^a_i$, from Equation~\eqref{eq:analysed_field_F}, using the particle set $\mathcal P_i^f$ directly.
\end{itemize}
The choice of the filter depends on the context, particularly regarding the feasibility of a remeshing process.

\subsection{Remesh-EnKF}\label{sec:remesh_enkf}
The first method consists of defining a scheme that combines a projection on a grid with a remeshing process to generate a new set of particles to discretize the assimilated states. Following Section~\ref{remesh_part}, the global scheme is built to conserve as much of the property of the forward members as possible. The assimilation is performed with the following steps:
\begin{enumerate}
    \item \textbf{Propagation}: Each member $i$ is independently updated thanks to the model, giving the new particle set $\mathcal{P}^f_i = \{(\bx^f_{i,p}, U^f_{i,p})\}_{p = 1}^{N_{i}}$,
    \item \textbf{Projection}: The associate field is projected on a uniform grid of characteristic length $\ell = 2dp$. Using the assignment operator in Equation~\eqref{eq:nodal_values}. We obtain for each member $i$ the vector of nodal values $\vu^{f,g}_{i}$.

    \item \textbf{Analysis}: Based on this new representation of the states, an Eulerian-based data assimilation is applied on the nodal state values $ \vu^{f,g}_{i}$ such that the analysis state $\vu_{i}^{a,g}$ is
          \begin{equation*}
              \vu^{a,g}_{i} = \vu^{f,g}_{i} + \sum_{j=1}^{N} F_{ji} \vu^{f,g}_{j},
          \end{equation*}
          The matrix $\bm F$ is given by \eqref{eq:analysed_field_F}.
    \item \textbf{Interpolation}: The analyzed grid values lead to the definition of the grid fields $u_i^{a,g}$ as in Equation~\eqref{eq:grid_approx}. A new uniform particle discretization is initialized. Two particles are placed (in 1D) inside each cell of the grid. The new particle intensities are approximate by Equation~\eqref{eq:new_set_intensities}, such that for $p \in \mathcal P^a$ the particle intensity is set
          \begin{equation*}
              U_{i,p}^a = u_i^{a,g}(\bx_p^a)~V_{p}.
          \end{equation*}
\end{enumerate}

\begin{figure}[htbp]
    \centering
    \makebox[\textwidth][c]{%
        \begin{subfigure}{0.40\textwidth}
            \centering
            \includegraphics[width=1.03\textwidth]{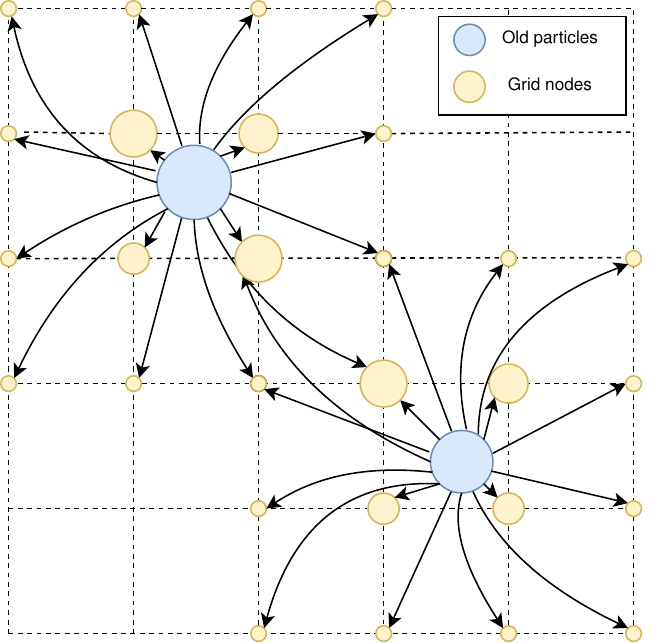}
            \caption{Projection}~\label{fig:remesh_enkf_1}
        \end{subfigure}%
        \begin{subfigure}{0.40\textwidth}
            \centering
            \includegraphics[width=\textwidth]{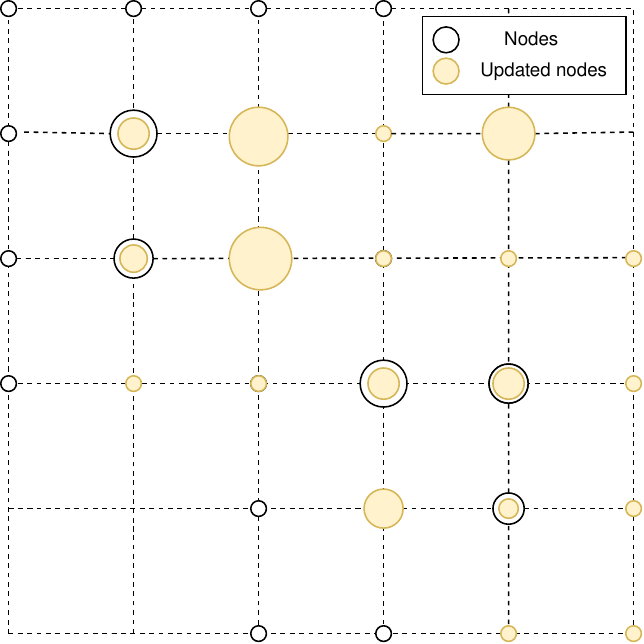}
            \caption{Analysis}~\label{fig:remesh_enkf_2}
        \end{subfigure}%
        \begin{subfigure}{0.40\textwidth}
            \centering
            \includegraphics[width=\textwidth]{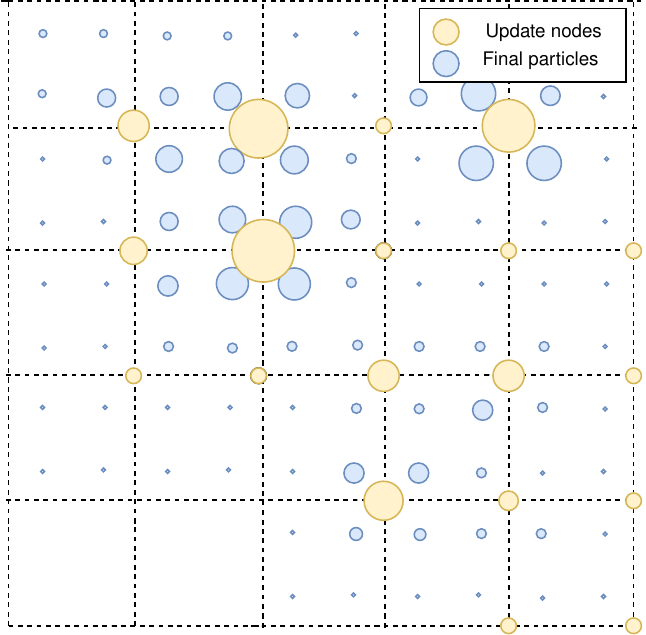}
            \caption{Interpolation}~\label{fig:remesh_enkf_3}
        \end{subfigure}}
    \caption{Illustration of the Remesh-EnKF. (\ref{fig:remesh_enkf_1}) The particles of the members are projected on a fixed Eulerian grid. (\ref{fig:remesh_enkf_2}) EnKF matrix operations are applied to nodal values. (\ref{fig:remesh_enkf_3}) New uniform sets of particles are generated with thresholding to interpolate the analyzed states.}~\label{fig:remesh_enkf}
\end{figure}

The Remesh-Filter update scheme is illustrated in Figure~\ref{fig:remesh_enkf}. We observe that all the operators involved in the update are linear with respect to the intensities of the particles and the nodal values. Therefore, performing the update on the grid or directly on the new particles is equivalent. However, the current algorithm allows an analysis on a lower-dimensional space, namely on the projection grid. It is also more practical for thresholding $\mathcal P^a$ along step 4. Specifically, the new particle is created if $|u_i^{a,g}(\bm x_p^a)| > \varepsilon_{\text{cut}}$, where $\varepsilon_{\text{cut}}>0$ the threshold value.

\subsection{Part-EnKF}\label{sec:method_part_enkf}
The previous adaptation of the EnKF filter allows the number of particles to be controlled during the analysis by regenerating a new particle distribution. In a way, it involves projecting the solution onto an Eulerian discretization before traditionally applying the analysis, followed by the reconstruction of the particle approximation. The Part-EnKF formulation aims to bypass the remeshing step, as this operation may not be feasible in all Lagrangian models. Specifically, the particle set of each member is kept the same during the update of the intensities. The fields, defined in~\eqref{eq:analysed_field_F}, are approximated using their respective particle sets $\mathcal P_i$. The Part-EnKF filter is built upon the following 3 steps:
\begin{enumerate}
    \item \textbf{Propagation}: Each member $i$ is independently updated thanks to the model, given the new particle set $\mathcal{P}^f_i = \{(\bx^f_{i,p}, U^f_{i,p})\}_{p = 1}^{N_{i}}$,
    \item \textbf{Analysis}: After the acquisition of a new observation and based on the current predictions, the analyzed fields are determined based on the discretization-free formula \eqref{eq:analysed_field_F} given for any coordinate $\bm x \in \Omega$
          \begin{equation*}
              u^a_{i}(\bm x) = u^f_{i}(\bm x) + \sum_{j=1}^{N} F_{ji} u^f_{j}(\bm x), \quad \text{where } u_i^f(\bx) = \sum_{p \in \mathcal P_i^f} U^f_{i,p} \phi_\varepsilon(\bm x - \bx_{i,p}^f)
          \end{equation*}
    \item \textbf{Approximation}: Based on the approximation formula we determined the new intensities $U^a_{i,p}$. This way, the analyzed field $u^a_i$ is approximated by
          \begin{equation*}
              u^a_i(\bm x) \simeq \sum_{p \in \mathcal P_i^f} U^a_{i,p} \phi_\varepsilon(\bm x - \bx^f_{i,p}), \quad \text{where } U^a_{i,p} = u^a_i(\bx_{i,p}^f)~V_p.
          \end{equation*}
          A more advanced approximation of $u^a_i$ can be used for the definition of $U^a_{i,p}$ based on regression operations discussed in Section~\ref{interpOp}.
\end{enumerate}

\begin{figure}[htbp]
    \centering
    \makebox[\textwidth][c]{%
        \begin{subfigure}{0.38\textwidth}
            \centering
            \includegraphics[width=\textwidth]{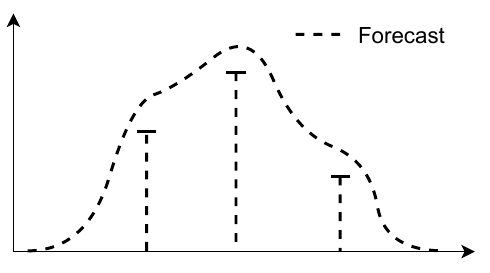}
            \caption{Forecast function}~\label{fig:part_enkf_1}
        \end{subfigure}%
        \begin{subfigure}{0.38\textwidth}
            \centering
            \includegraphics[width=\textwidth]{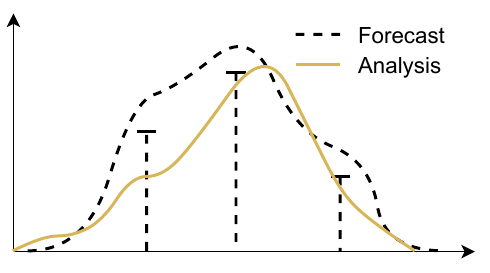}
            \caption{Analysed function}~\label{fig:part_enkf_2}
        \end{subfigure}%
        \begin{subfigure}{0.38\textwidth}
            \centering
            \includegraphics[width=\textwidth]{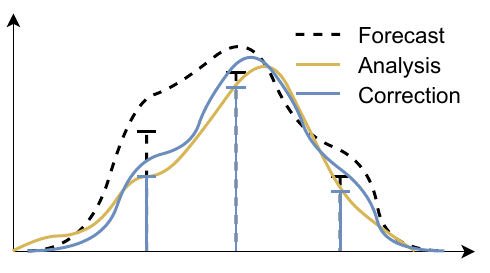}
            \caption{Intensity correction}~\label{fig:part_enkf_3}
        \end{subfigure}}%
    \caption{One-dimensional illustration of the Part-EnKF method. (\ref{fig:part_enkf_1}) A set of three particles discretizes a member forecast function. (\ref{fig:remesh_enkf_2}) The analyzed solution is determined thanks to equation~\eqref{eq:analysed_field_F}. (\ref{fig:part_enkf_3}) The forecast particle intensities are updated to fit the analyzed solution.}~\label{fig:part_enkf}
\end{figure}

The Part-EnKF filter is illustrated in Figure~\ref{fig:part_enkf}. The complexity is in $O(N^2)$ when the members have comparable sizes of particle sets. This complexity can be considerably reduced when working with compact $\phi_\varepsilon$, and relying on fast particle algorithms, such as fast Gaussian summation and multipole expansion~\cite{GREENGARD1987325,lee_faster_2012}.

\section{1D density advection-diffusion problem}\label{App_1D}
\subsection{Problem setting}
An initial exploration is conducted on a one-dimensional application to assess the performance of the filters. We define the following one-dimensional $2\pi$-periodic advection-diffusion problem; consider
\begin{equation}\label{eq:adv_diff}
    \frac{\partial u}{\partial t}(x,t) + v \frac{\partial u}{\partial x}(x,t) = \visc \frac{\partial^2 u}{\partial x^2}(x,t),
\end{equation}
with $x$ the spatial coordinate, $v$ a constant advection velocity and $\visc >0$ a constant diffusion coefficient.
In the following, the reference solution will use $v = \refv$ and $\visc = \refvisc$ as parameters.
We define the $2\pi$-periodic heat kernel in one dimension, such that
\begin{equation*}
    K(x, t) = \sum_{k=-\infty}^{\infty} \frac{1}{\sqrt{4 \pi t}} \exp{\left(-\frac{{(x - 2\pi k)}^2}{4t} \right)}.
\end{equation*}

\begin{figure}[htbp]
    \centering
    \includegraphics[width=\textwidth]{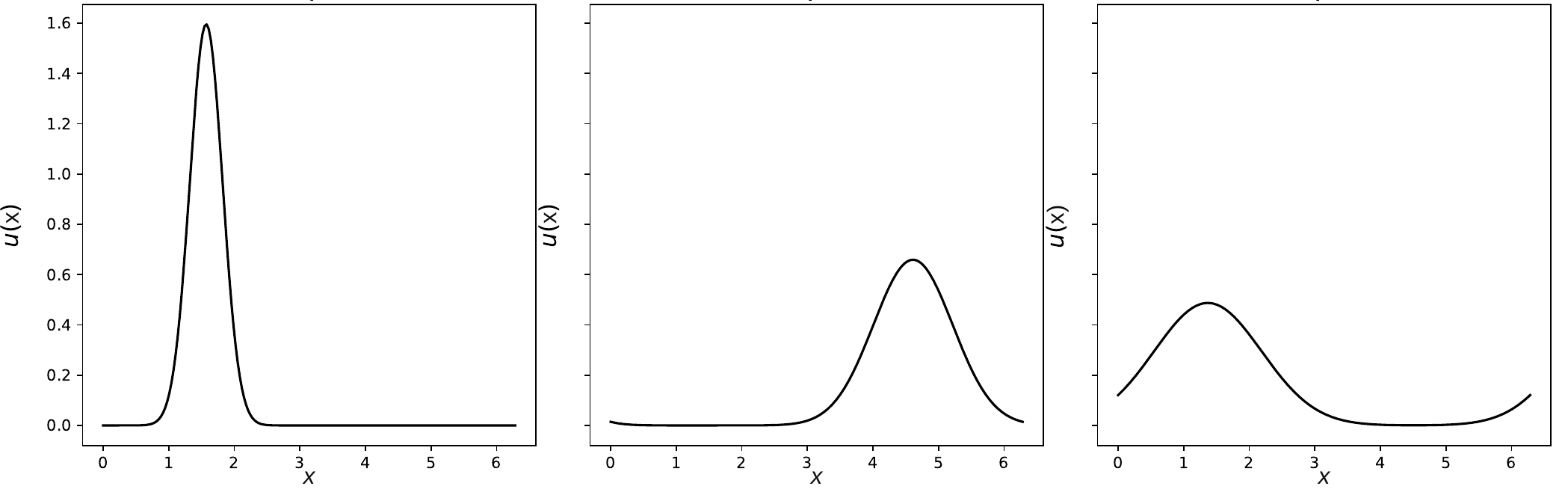}
    \caption{The exact solution of the advection-diffusion problem at different times: $t=0$, $t= \frac{\pi}{v}$, and $t= \frac{2\pi}{v}$.}
    \label{fig:1d_analytical}
\end{figure}

Considering an initial condition defined by $u^{gt}(x, 0) = K(x-x_0, \sigma_0^2 / 2)$, where $\zz$, and $\sigz$, the exact analytical solution is
\begin{equation*}
    u^{gt}(x, t) = K(x - v t - x_0, \visc t + \sigma_0^2 / 2).
\end{equation*}The analytical solution has succinctly a Gaussian shape centered on a point that moves with advection velocity and a width increasing in time as $\sqrt{Dt}$. Figure~\ref{fig:1d_analytical} shows this solution at different times.

The Lagrangian discretization of~\eqref{eq:adv_diff} with particles $\left\{x_p(t), U_p(t)\right\}_{p=1}^{N_p}$ yields to governing equations
\begin{equation*}
    \frac{dx_p(t)}{dt} = v(x_p(t), t) = v, \quad \frac{dU_p(t)}{dt} \simeq V_p \left(D \frac{d^2 u(x_p(t))}{dx^2}\right)
\end{equation*}

To solve the advection-diffusion equation, we employ the two steps of the viscous splitting algorithm. First, the advection is taken into account by moving the particle's position with the velocity $\bm v$.
Second, we use a redistribution method called the Particle Strength Exchange Method (PSE)~\cite{degond_1989,cottet_1990} to account for the diffusion. It comes
\begin{equation*}
    \frac{dU_p(t)}{dt} = D \varepsilon^{-2} \sum_{q \in \mathcal P} (V_p U_q(t)-V_q U_p(t)) \overline{\overline{\phi}}_\varepsilon( x_q(t) - x_p(t)),
\end{equation*}where $V_p$ is the volume of the particle $p$. The diffusion kernel $\overline{\overline{\phi}}_\varepsilon$ is here taken to be the square exponential kernel. For further details on the computation, please refer to~\cite{cottet_1990}.

For the periodic boundary problem described in~\eqref{eq:adv_diff}, we define an equivalent kernel function
\begin{equation*}
    \phi_\varepsilon(x) = \sum_{n=-\infty}^{+\infty} \phi^g_\varepsilon(x - 2 \pi n ), \quad \text{given } u(x) = \sum_{p \in \mathcal P} U_p(t) \phi_\varepsilon(x - x_p(t)),
\end{equation*}where $\phi^g_\varepsilon$ is the squared exponential kernel on an infinite domain.

\subsubsection{Assimilation settings}
\begin{figure}[htbp]
    \centering
    \begin{subfigure}[t]{0.49\textwidth}
        \includegraphics[width=\textwidth]{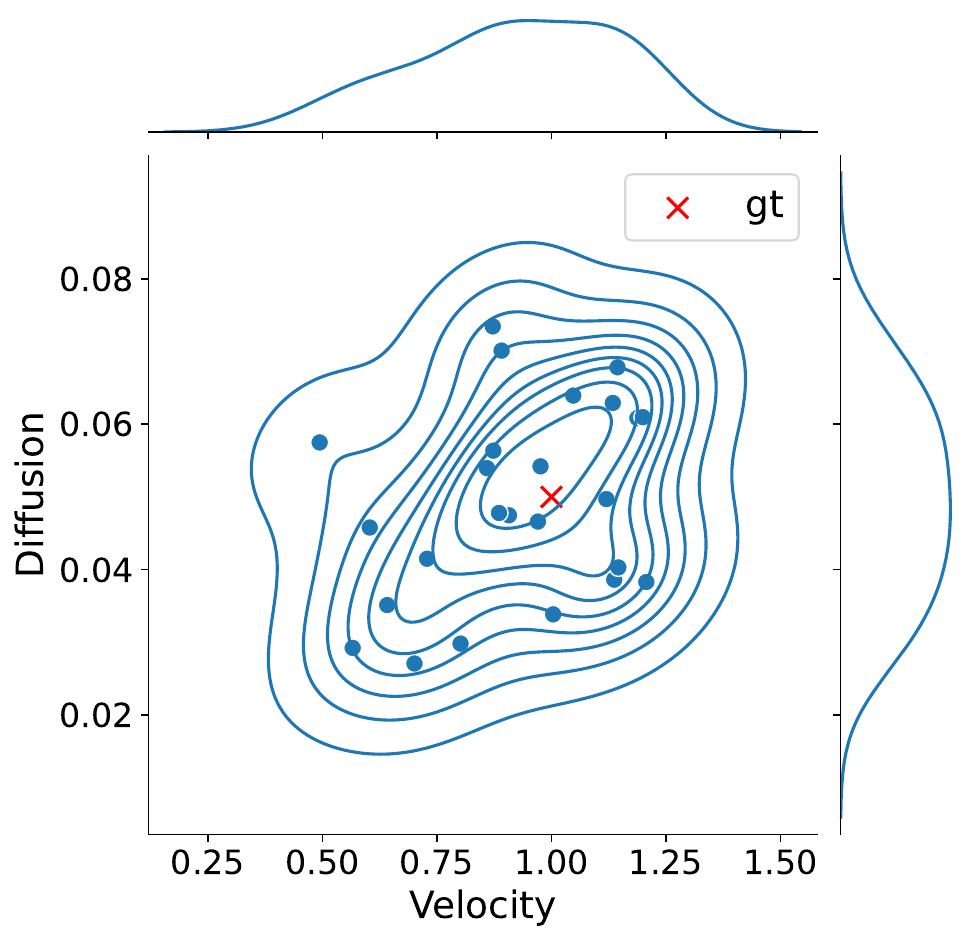}
    \end{subfigure}%
    \hfill
    \begin{subfigure}[t]{0.49\textwidth}
        \includegraphics[width=\textwidth]{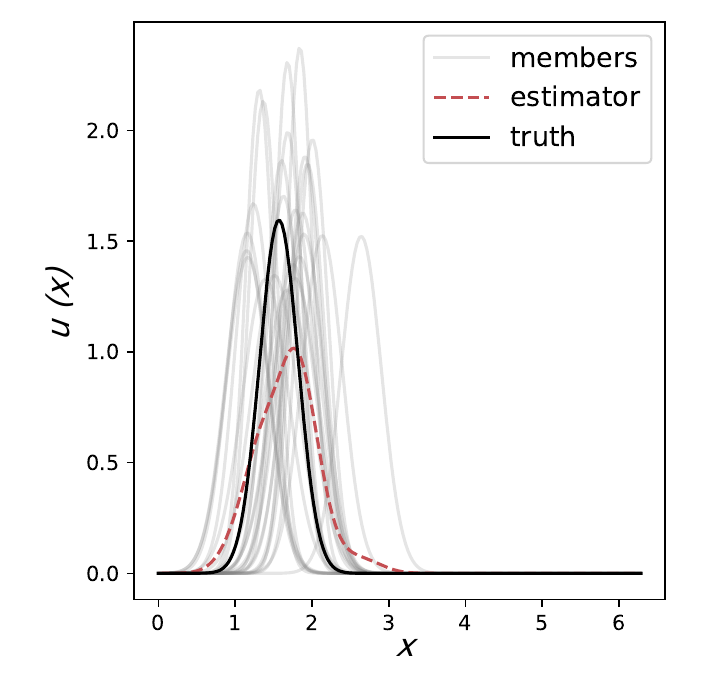}
    \end{subfigure}%
    \caption{On the left, the initial parameter sample with $v$ on the abscissa and $D$ on the ordinate. On the right, the initial ensemble states.}\label{fig:initial_gen}
\end{figure}

All filters are tested on an identical initial prior ensemble
\begin{equation*}
    \left\{u_i(x, t=0) = K(x - x_0^{(i)}, (\sigma_0^{(i)})^2 / 2) \right\}_{i=1}^N
\end{equation*}with $N = 25$, where $x_0^{(i)}$ et $\sigma_0^{(i)}$ are randomly sampled from their distributions reported in Table~\ref{tab:ens_gen_1d}. The velocity $v^{(i)}$ and diffusion $D^{(i)}$ are also sampled for the considered model. They will be integrated into the data assimilation scheme to perform calibration of model parameters. The initial parameters samples and initial state are illustrated in Figure~\ref{fig:initial_gen}.
\begin{table}[htbp]
    \centering
    \caption{Ensemble generation variables}
    \begin{tabular}[t]{|l|l|}
        \hline
        Variables & Distributions                                   \\
        \hline
        position  & $x_0^{(i)} \sim \mathcal{N}(\meanZm, \sigmaZm)$ \\
        spread    & $\sigma_0^{(i)} \sim\mathcal{U}(\smLow, \smUp)$ \\
        velocity  & $v^{(i)} \sim \mathcal{N}(\vmean, \vstd)$       \\
        diffusion & $D^{(i)} \sim \mathcal{U}(\Dlow, \Dup)$         \\
        \hline
    \end{tabular}
    \label{tab:ens_gen_1d}
\end{table}

We conduct $N_{\text{assim}} = 30$ assimilation steps at evenly spaced intervals until the final time $t_f = \frac{2 \pi}{v}$. In each assimilation step, the ground truth $u^{gt}$ is observed at 6 evenly spaced positions. The observations' noise is assumed Gaussian, independent, with a fixed variance $\sigma_\varepsilon^2 = 0.05$.

\subsubsection{Numerical experiments}
We employ particles of size $d_p = \frac{2 \pi}{100}$ and a smoothing length of $\varepsilon = 1.3 d_p$.
In the particle-based simulation, the discretized uses regularly spaced particles. For this 1D problem, the advection does not involve distortion of the particle distribution. 
A parameter $\varepsilon_{\text{cut}}$ is introduced to remove particles with low intensities at initialization. This parameter is used to restrict the spatial extension of the particle set around $x_0^{(i)}$.

When using a purely Eulerian method, we use a finite difference solver discretized on a uniform grid using $N_{\text{grid}} =100$ nodes. The purely Eulerian filter is called Grid-EnKF. It consists of applying EnKF update to the nodal values.

\subsubsection*{Error definition}
We define the relative RMSE as the following relative ratio
\begin{equation*}
    \text{rRMSE} = \frac{1}{\norm{u^{gt}}_{L_2}} \left[\frac1N \sum_{i=1}^{N} \norm{u_i - u^{gt}}^2_{L_2}\right]^{1/2} \quad \text{where } \norm{u}^2_{L_2} = \int_0^{2\pi} |u|^2 dx.
\end{equation*}
The $L_2$ norm is computed using a numerical quadrature of a uniform fine grid over $[0, 2\pi]$.
Additionally, we compute the rRMSE for each parameter $v$ and $D$ as
\begin{equation*}
    rRMSE_{\theta}= \frac{1}{|\theta^{gt}|}\left[\frac1N \sum_{i=1}^{N} |\theta_i - \theta^{gt}|^2\right]^{1/2}.
\end{equation*}

\subsection{Results}
\begin{figure}[htbp]
    \centering
    \begin{subfigure}{\textwidth}
        \centering
        \includegraphics[width=0.8\textwidth]{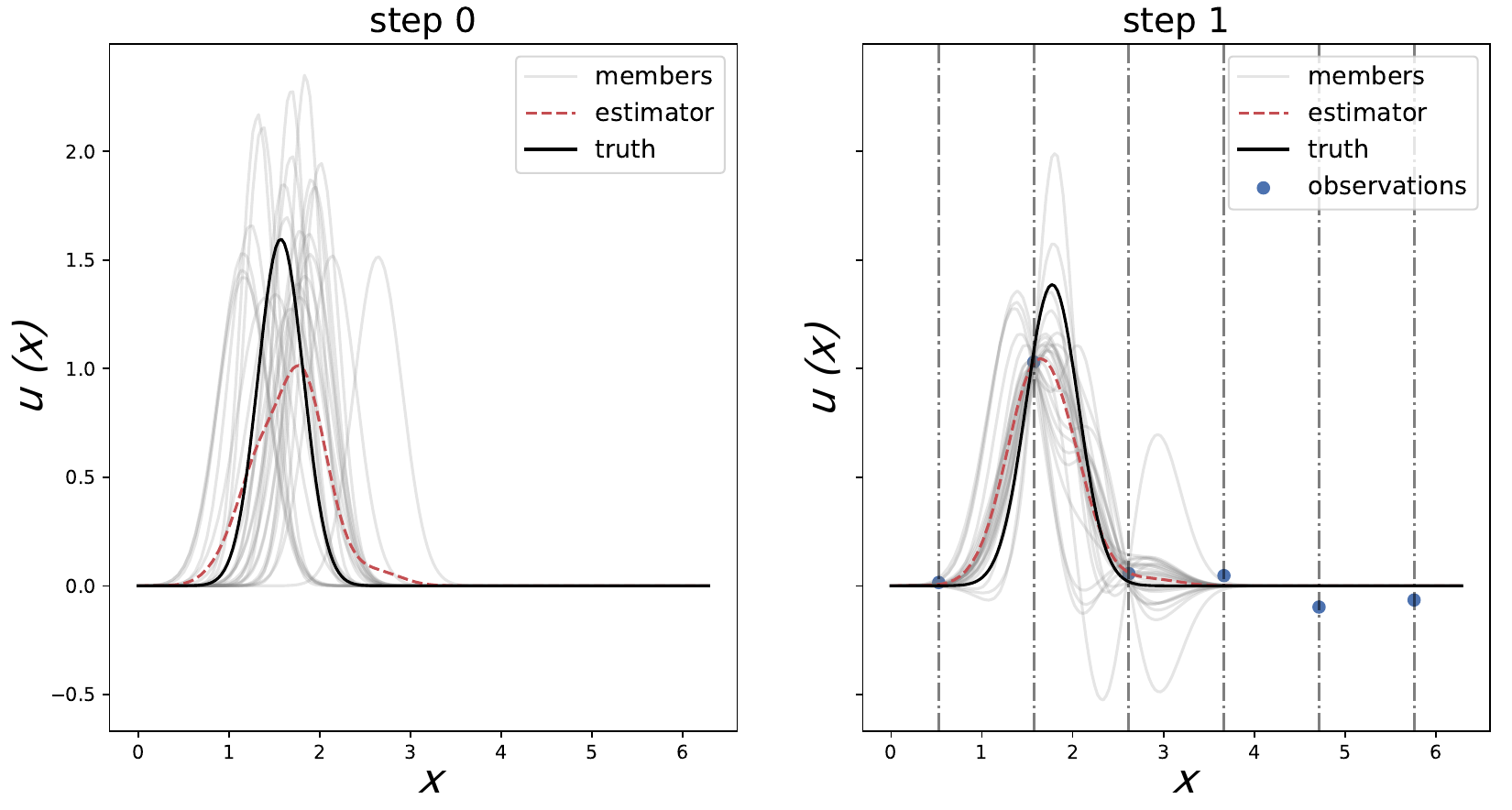}
    \end{subfigure} \\
    \vfill
    \begin{subfigure}{\textwidth}
        \centering
        \hspace{0.1cm}
        \includegraphics[width=0.8\textwidth]{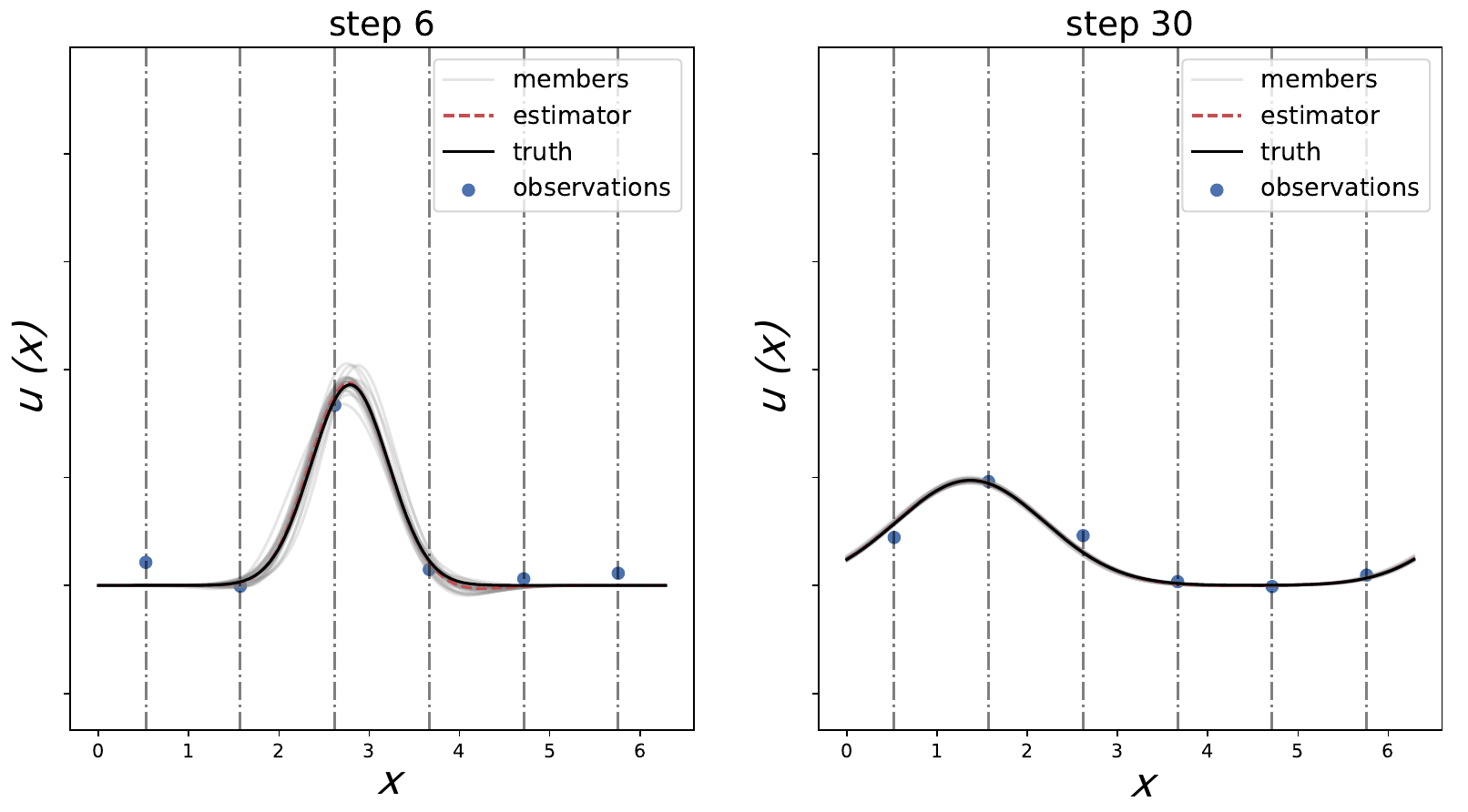}
    \end{subfigure}
    \caption{Analysed solution at different assimilation steps for the Remesh-EnKF filter.}\label{fig:remesh_assim_step}
\end{figure}

\begin{figure}[htbp]
    \centering
    \includegraphics[width=0.85\textwidth]{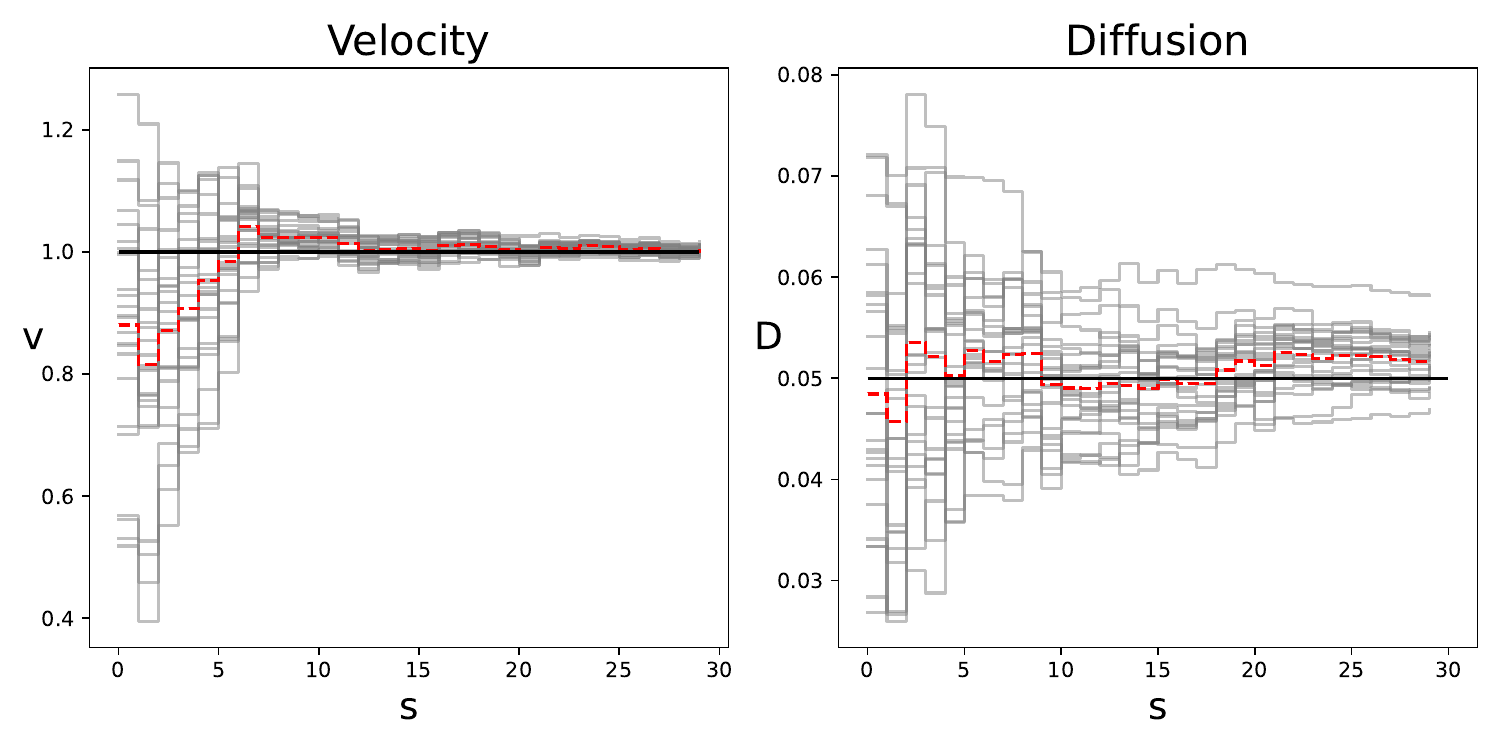}
    \caption{Calibration of model parameters over assimilation steps for the Remesh-EnKF filter.}\label{fig:parameter_time}
\end{figure}

We compare the performance of the different filters — Grid-EnKF, Remesh EnKF, and Part-EnKF — in state assimilation and parameter calibration (velocity and diffusion). Figure~\ref{fig:remesh_assim_step} shows the evolution of the solution of the Remesh-EnKF over time. The corrections at the first assimilation induce significant adjustments that alter the shape of the members with positivity loss. As the assimilation advances, the ensemble converges well toward the ground truth. Figure~\ref{fig:parameter_time} illustrates the evolution of the members' values of the two model parameters for the Remesh-EnKF. For each member, the parameters have piecewise constant solutions with jumps at the assimilation step. A reduction in the ensemble variance is clearly observed as well as a progressive convergence to the ground truth. The reduction in the variance of the advection speed is faster, as observations are more sensitive to it. The other filters have similar behavior.

\begin{figure}[htbp]
    \centering
    \includegraphics[width=0.8\textwidth]{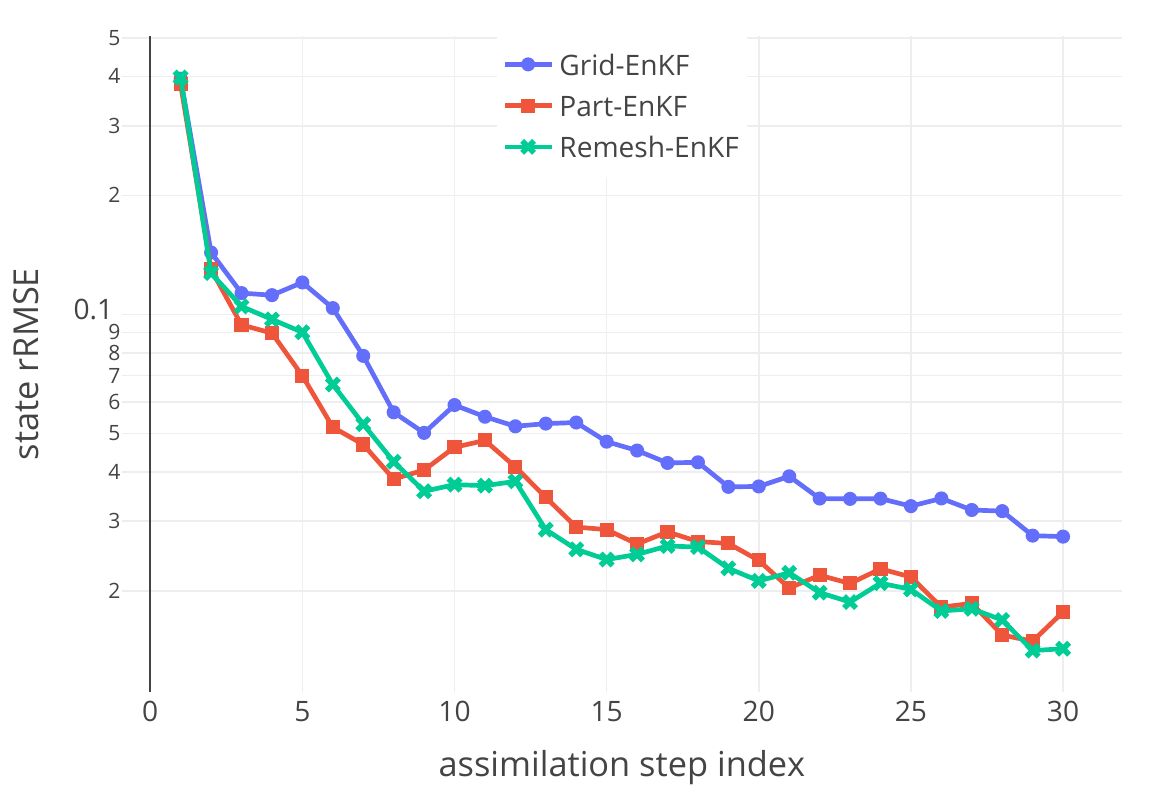}
    \caption{State rRMSE with respect to assimilation step index.}\label{fig:1d_error_time_state}
\end{figure}

To appreciate quantitatively the different filter results, Figure~\ref{fig:1d_error_time_state} shows the evolution of the rRMSE for the different methods. The errors behave similarly for all filters, with a fast reduction of the error during the first assimilation steps, followed by a significantly slower reduction rate until reaching a plateau. The two particle-based filters agree well with each other, with very similar error levels. The reference Eulerian filter offers slightly higher rRMSE, probably due to the solver scheme.

Concerning the model parameters, we also report similar reductions of rRMSE for the different filters. As already said, the rRMSE on velocity decreases rapidly during the first 6 assimilation steps, suggesting a strong reduction in bias and variances, while the rRMSE on diffusion decreases at a slower constant rate over the range of assimilation steps shown. Consequently, the two proposed filters offer an effective way to apply the EnKF correction scheme for particle-based discretizations in this 1D problem.
\begin{figure}[htbp]
    \centering
    \includegraphics[width=0.8\textwidth]{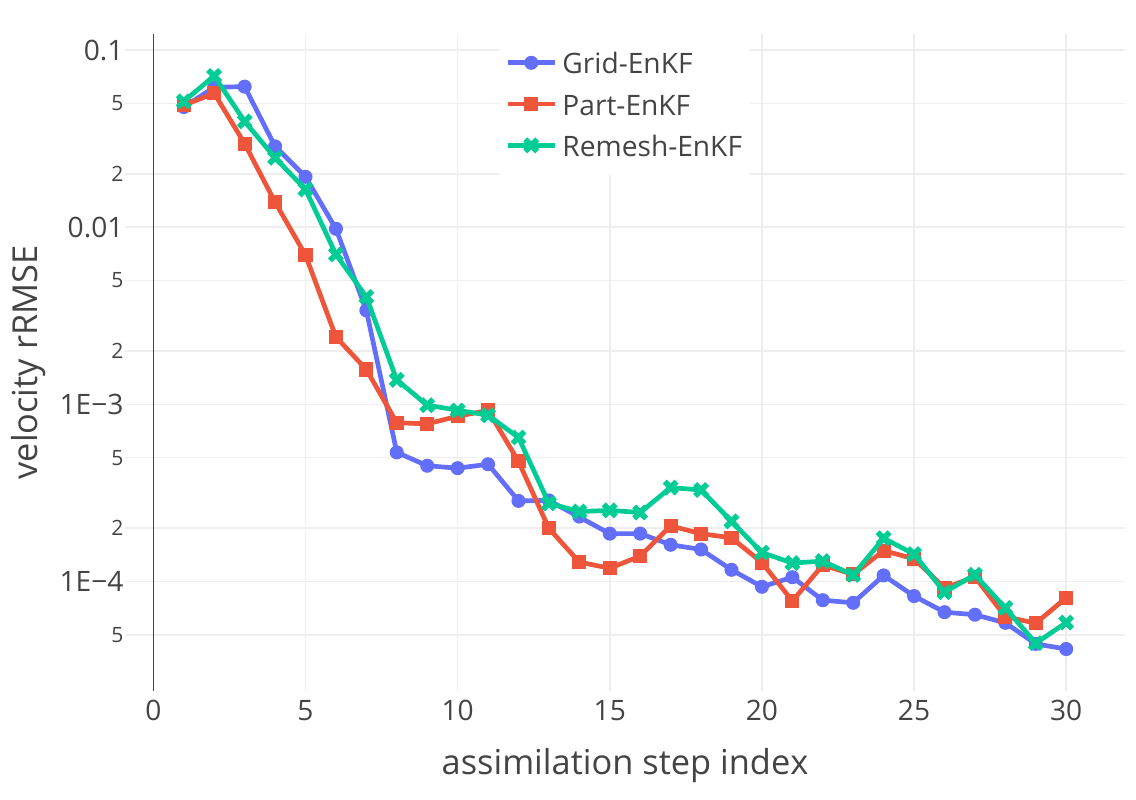}
    \caption{Velocity rRMSE with respect to assimilation step index.}\label{fig:1d_error_time_velocity}
\end{figure}%
\begin{figure}[htbp]
    \centering
    \includegraphics[width=0.8\textwidth]{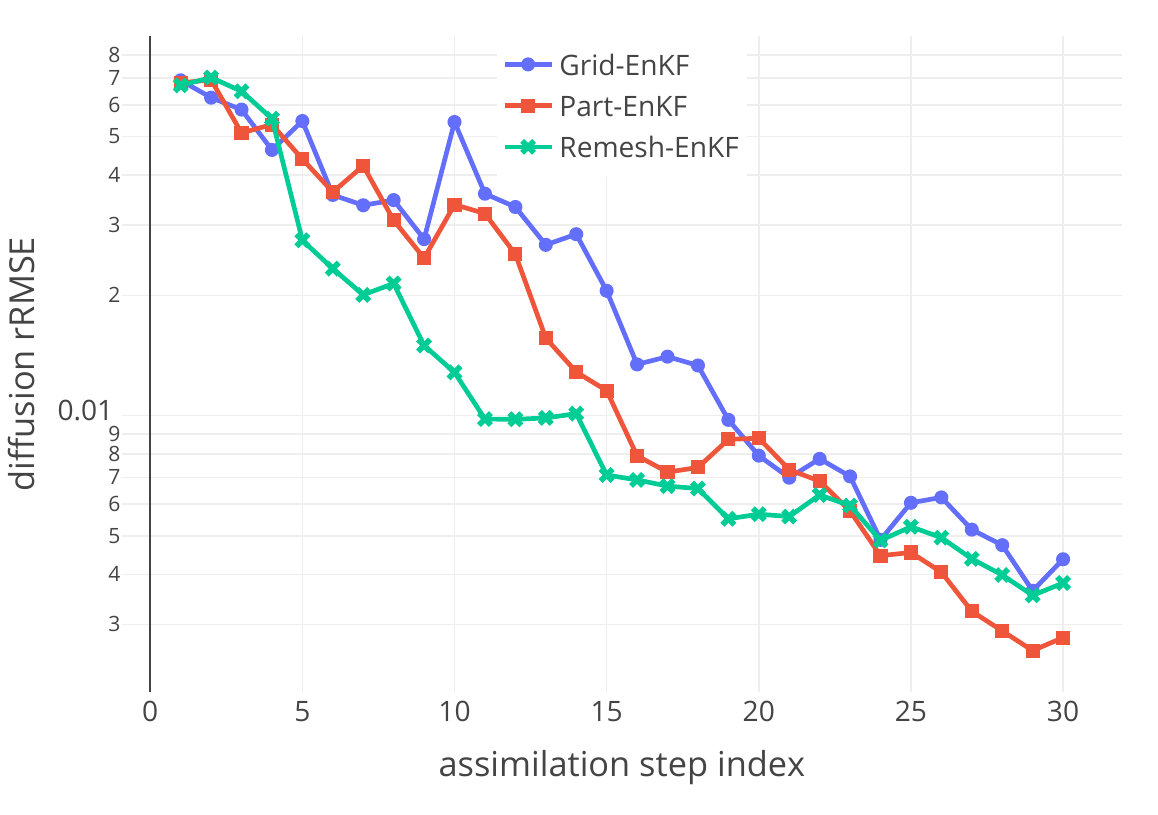}
    \caption{Diffusion rRMSE with respect to assimilation step index.}
    \label{fig:1d_error_time_viscosity}
\end{figure}%

\subsubsection{Influence of the support size}
By changing the thresholding $\varepsilon_{\text{cut}}$, we investigate the influence of the discretized support size on the filter solution. Without thresholding, the number of particles for each member is $N_{\text{part}}=100$. For $\varepsilon_{\text{cut}} = 1.0e-11$, the mean number of particles per member is $N_{\text{part}}=60$. Although the change in the particle approximation is negligible, the particle discretization covers only roughly 2/3 of the domain. By varying $\varepsilon_{\text{cut}}$, we control the extent of the particle discretization of the members.

\begin{figure}[htbp]
    \centering
    \includegraphics[width=0.8\textwidth]{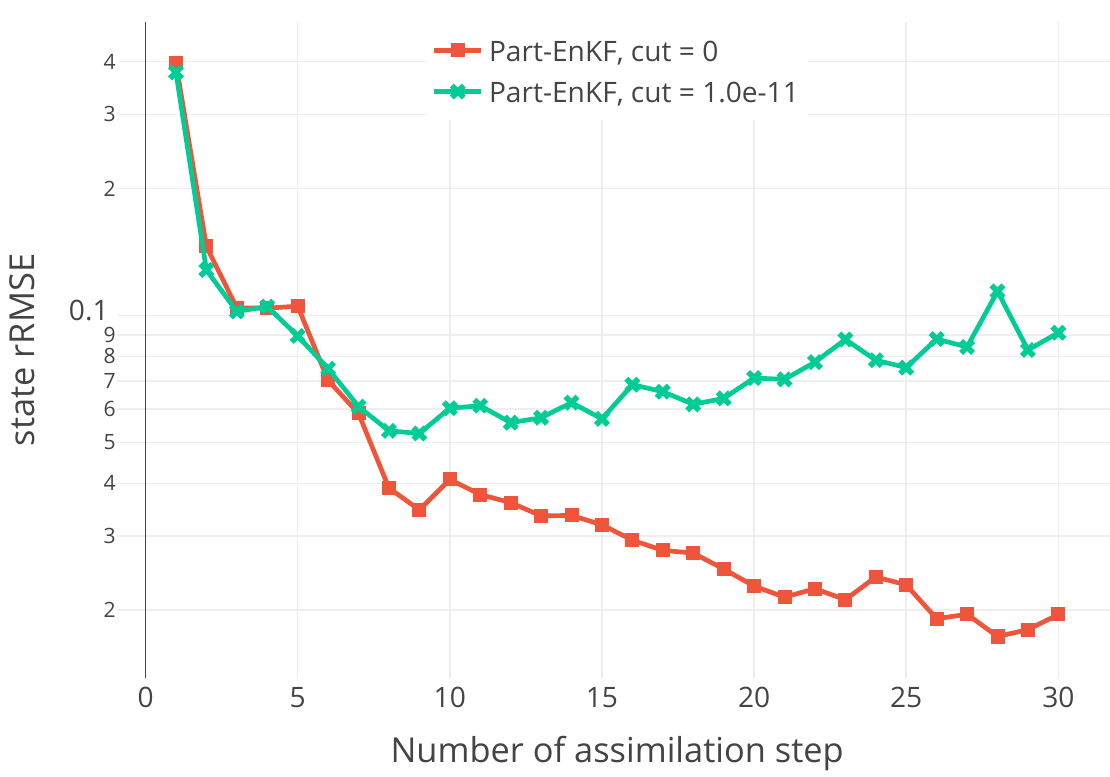}
    \caption{rRMSE on the state through time for the Part-EnKF method with different $\varepsilon_{\text{cut}}$ values.}\label{fig:part_extend}
\end{figure}
\begin{figure}[htbp]
    \centering
    \begin{subfigure}{0.85\textwidth}
        \centering
        \includegraphics[width=\textwidth]{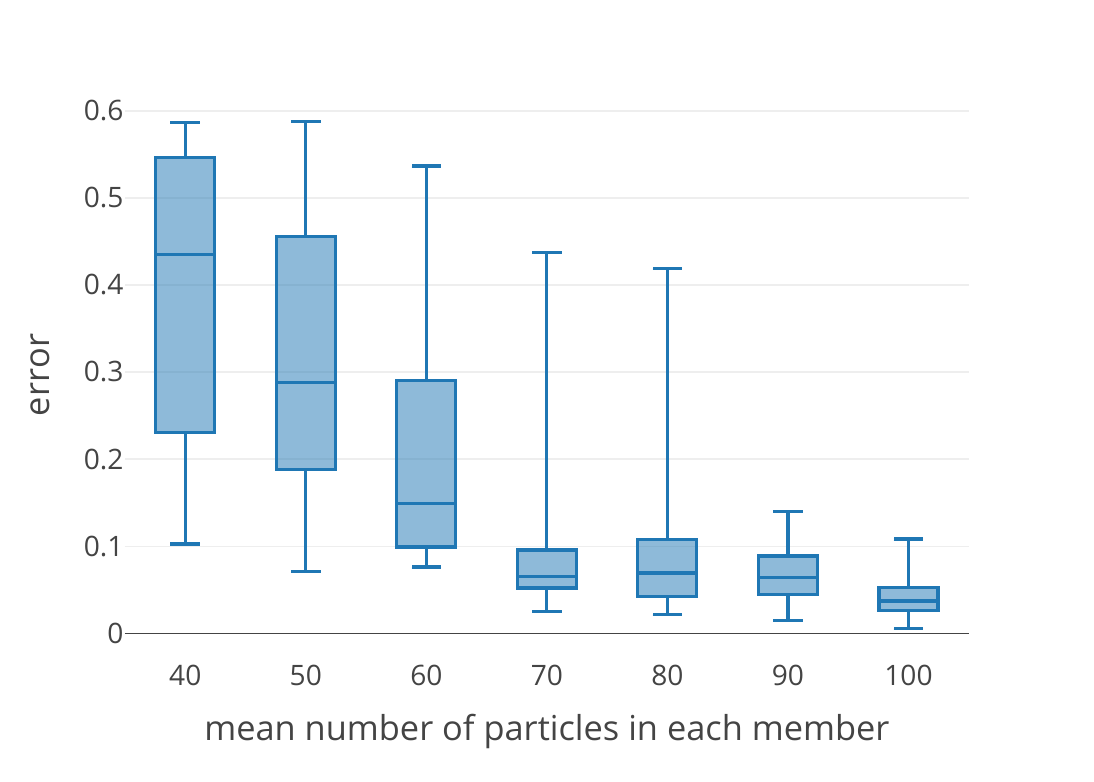}
    \end{subfigure}%
    \caption{Box plot of the normalized state discrepancy over the ensemble members for Part-EnKF. The threshold $\varepsilon_{\text{cut}}$ is varied to change the effective size of the discretization supports, here measured as the average number of particles in the members' discretization.}\label{fig:error_support1}
\end{figure}
\begin{figure}[htbp]
    \centering
    \begin{subfigure}[t]{0.49\textwidth}
        \centering
        \includegraphics[width=\textwidth]{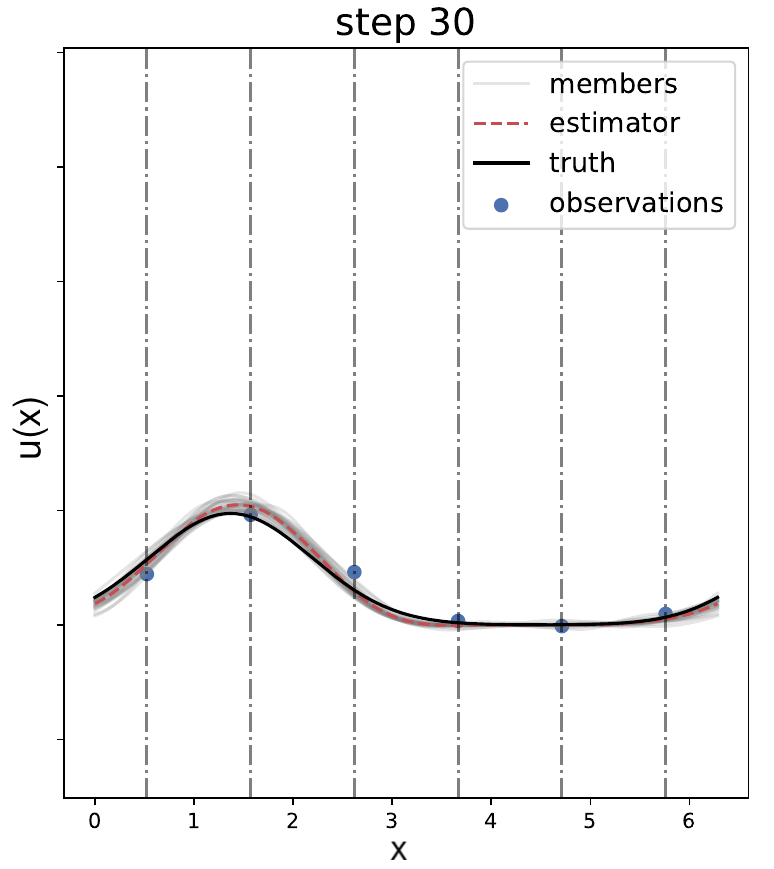}
        \label{error_support2}
    \end{subfigure}%
    \begin{subfigure}[t]{0.49\textwidth}
        \centering
        \includegraphics[width=\textwidth]{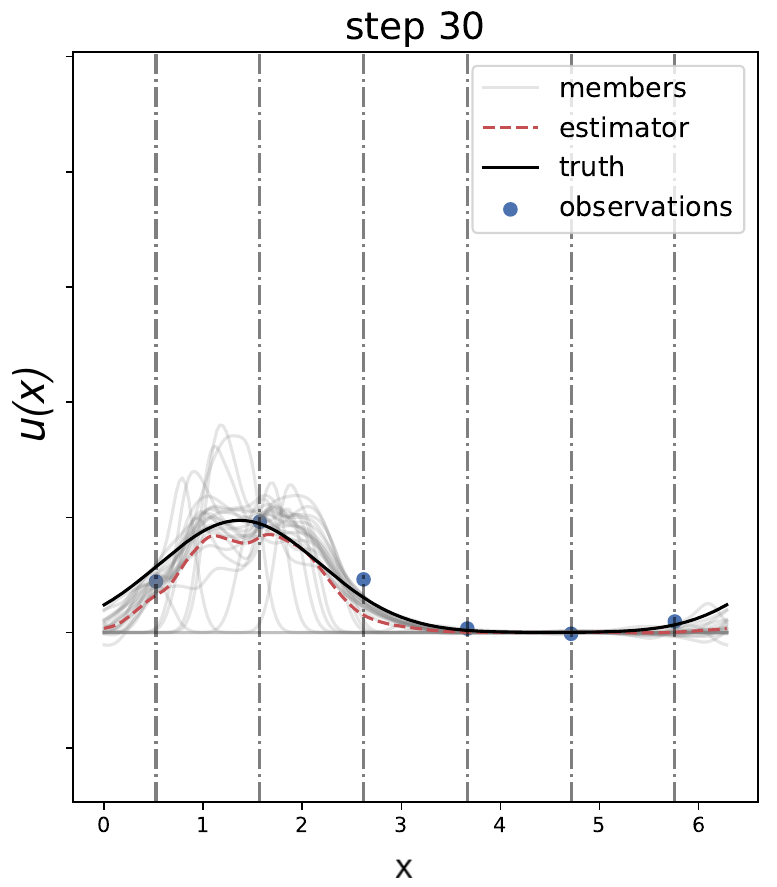}
        \label{error_support3}
    \end{subfigure}
    \caption{State ensemble at final assimilation step. Left: $\varepsilon_{\text{cut}} = 0$ (100 particles), Right: $\varepsilon_{\text{cut}} = 1.0e-11$ (60 particles).}
    \label{fig:error_support}
\end{figure}

In Figure~\ref{fig:part_extend}, we observe that Part-EnKF with a more restricted discretization is not able to converge adequately. Increasing the number of assimilations does not change this behavior. The origin of the error for $\varepsilon_{\text{cut}}$ is the quadrature error appearing when defining the analyzed solution on the members' discretization. The restricted support of the member discretization does not sufficiently overlap with the support of the analyzed solution, inducing errors. Even with adequate regression, the projection of the analyzed solution remains based on the forecast support. It is imperative to increase the discretized support of the members' particles to achieve a better approximation.
We validate this explanation by varying the size of the discretized support, which is reported here as the average number of particles per member. Quantitatively, as observed in Figure~\ref{fig:error_support1}, the error estimate and intra-ensemble dispersion decrease in the number of particles involved in the approximation increases. At the level of 70 particles discretization (i.e., 70\% of the domain is covered by particles), a threshold is reached. Moreover, quantitatively examining Figure~\ref{fig:error_support}, the right plot reveals that the solution closely aligns with the reference but is unable to approximate the solution appropriately when the discretization support is not as large. On the other hand, the calibration of the parameters is still adequately performed, as the dynamics of the problem are well represented (i.e., the predictions of the center and the width of the field solution). While the calibration error is not explicitly shown here, it follows the same patterns as shown in Figure~\ref{fig:parameter_time}.

As already mentioned, the prior ensemble was adequately represented by the restricted particle support. Adding particles along the assimilation is the only way to have a better solution approximation. An alternative would be to change the position of the current particles.
This example underscores the ability of both filters to yield results comparable, or even better, to a classical purely Eulerian Grid-EnKF filter. Additionally, we highlight the Part-EnKF capability of assimilating on a given particle discretization while also emphasizing the importance of addressing discrepancies between members, which can pose challenges to approximate analyzed solutions.
\newpage
\section{Application to incompressible flow via Vortex Method}\label{App_2D}
\subsection{Problem setting}
\subsubsection{Model equation}
In this section, we apply the Vortex Method to a two-dimensional scenario, as outlined by~\cite{cottet_vortex_2000}. The Vortex Method is a Lagrangian approach utilizing particles to discretize the vorticity field, allowing for the solution of the Navier-Stokes equation for viscous incompressible flow where $\omega~\bm z = \nabla \times \bm{v}$ satisfies
\begin{equation}\label{eq:vorticity_flow}
    \frac{\partial \omega}{\partial t} + (\bm{v} \cdot \nabla) \omega - \nu \Delta \omega = 0,
\end{equation}where $\omega = \frac{\partial v_y}{\partial x} - \frac{\partial v_x}{\partial y}$ denotes the scalar vorticity field normal to the flow plan, $\bm{v}$ is the velocity, and $\nu>0$ is the viscosity of the fluid. The vorticity field is discretized using a collection of discrete vortices, each characterized by a position $\bm x_p$, a volume $V_p$, an associated kernel $\phi_\varepsilon$, and a circulation $\Gamma_p$. For all points $\bm x$ within the computational domain $\Omega \subseteq \mathbb R^2$, the vorticity is expressed as
\begin{equation*}
    \omega(\bm x,t) = \sum_{i=1}^{N_p} \Gamma_p(t) \phi_\varepsilon(\bm x - \bm x_p(t)).
\end{equation*}

To solve the convection-diffusion equation, we employ a viscous splitting scheme, following the methodology outlined in~\cite{cottet_1990}, acknowledging the predominance of the convection term over viscosity. We use the Vortex-In-Cell algorithm~\cite{christiansen_1973, birdsall_1969}, coupled with an FFT solver to compute the advection velocity from the stream-function $\psi$ verifying $\Delta \psi = - \omega$. The algorithm involves assigning vorticity values on an underlying uniform grid using a particle-to-grid formula in order to solve the Poisson equation for the stream function on the grid. After the differentiation of the stream function to get the velocity components at the grid nodes, interpolation schemes are employed to compute the velocity at the location of the particles. A Runge-Kutta 3 time-stepping scheme is employed to integrate the particle positions transported with this velocity field. The diffusion step is solved with the PSE methods previously described in Section~\ref{App_1D}.
Using the Lagrangian formulation of the equation~\ref{eq:vorticity_flow}, the dynamic equation for the particles becomes
\begin{gather*}
    \left\{\begin{aligned}
         & \frac{d \bx_p(t)}{dt} = \bm v(\bx_p(t), t),                                                                                                       \\
         & \frac{d\Gamma_p(t)}{dt} = \nu \varepsilon^{-2} \sum_{q \in \mathcal P} [V_p \Gamma_q(t) - V_q \Gamma_p(t)] \phi_\varepsilon(\bx_p(t) - \bx_q(t)).
    \end{aligned}\right.
\end{gather*}

\subsubsection{Lamb-Chaplygin dipole and simulation parameters}
\begin{figure}[htbp]
    \centering
    \includegraphics[width=0.65\textwidth]{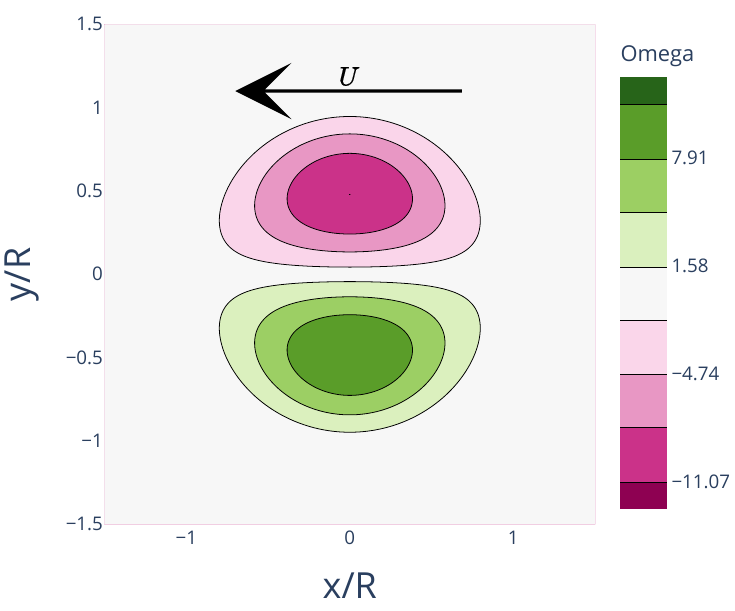}
    \caption{The Lamb-Chaplygin dipole vorticity field on a normalized space.}
    \label{fig:lamb_dipole}
\end{figure}

We define the reference state as the advection of the Lamb-Chaplygin dipole inside a close domain with stress-free walls. The Lamb-Chaplygin dipole is a popular choice for numerical studies~\cite{orlandi_vortex_1990}. The dipole is characterized by a translation velocity $U$, a mean position $\bm{x}_0$, a radius $R$, and an orientation $\alpha$. The dipole vorticity field $\omega$ is expressed as
\begin{equation*}
    \omega(r) = \begin{cases}
        \frac{-2 k U J_1(kr)}{J_0(kR)} \sin \alpha \quad & \text{for} \quad r < R, \\
        0 \quad                                          & \text{otherwise},
    \end{cases}
\end{equation*}where $(r, \alpha)$ are the polar coordinates with respect to the dipole's center. Here, $J_0$ and $J_1$ denote the zeroth and first-order Bessel functions of the first kind, respectively, and $k$ is determined such that $kR$ corresponds to the first non-trivial zero of the first Bessel function. The dipole vorticity field is depicted in Figure~\ref{fig:lamb_dipole}. Here, the evolution is more complex because of the presence of the boundaries and the diffusion.

\begin{figure}[htbp]
    \begin{subfigure}[b]{0.33\textwidth}
        \includegraphics[width=\textwidth]{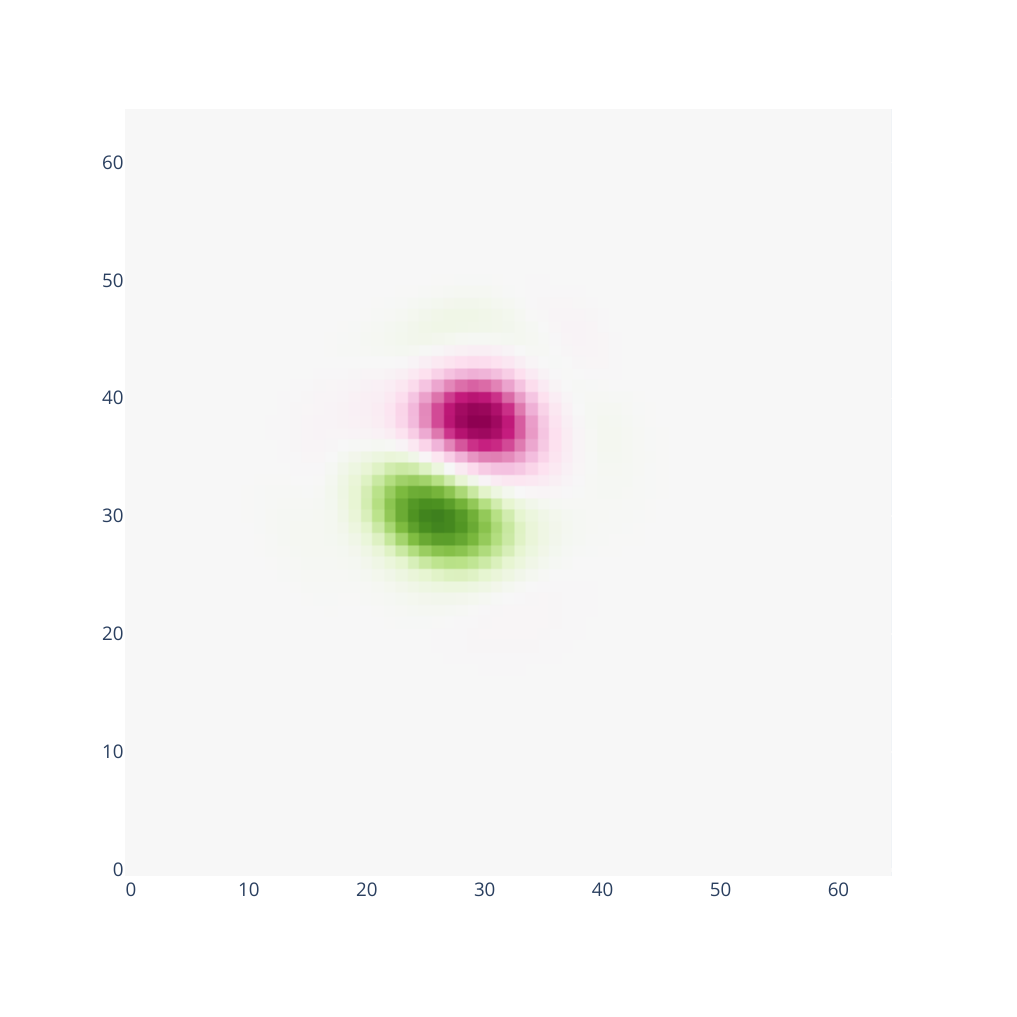}
    \end{subfigure}%
    \hfill
    \begin{subfigure}[b]{0.33\textwidth}
        \includegraphics[width=\textwidth]{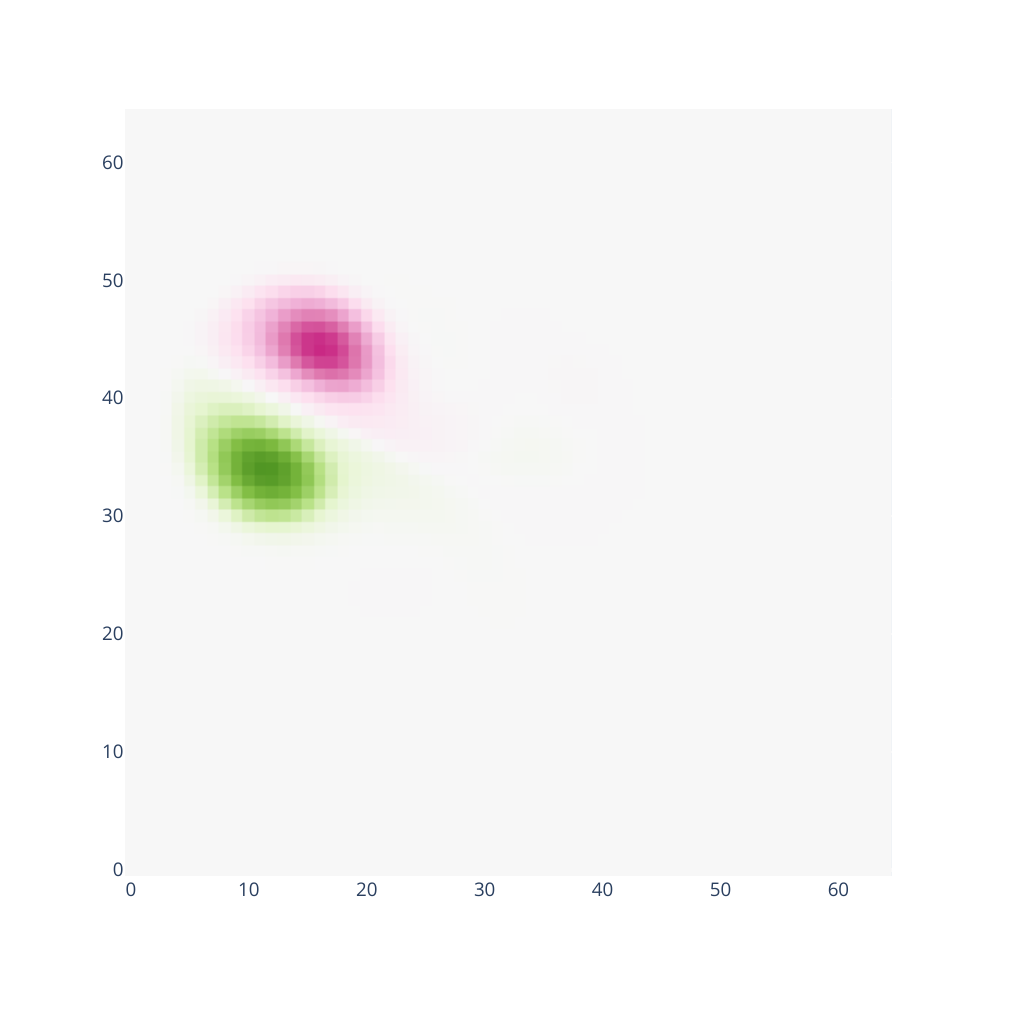}
    \end{subfigure}%
    \hfill
    \begin{subfigure}[b]{0.33\textwidth}
        \includegraphics[width=\textwidth]{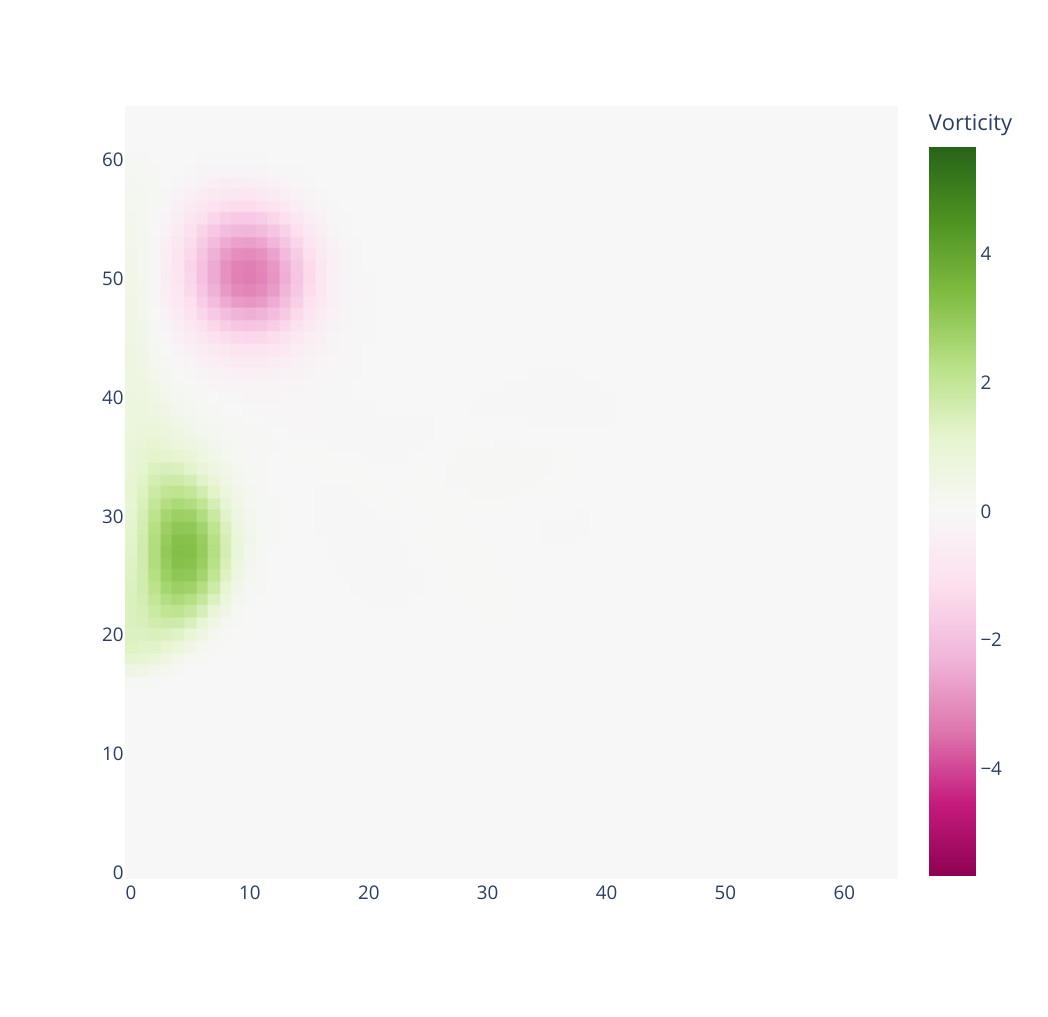}
    \end{subfigure}%
    \caption{Evolution of the ground truth. Vorticity field at $t=1, 5$ and $10$.}
    \label{fig:ref_trajectory}
\end{figure}

The dipole is positioned at the center of a box with dimensions $[0, \pi]^2$, featuring an orientation of $\alpha = \frac{7\pi}{8}$ rad., a radius of $R = 0.5$ meters, and a velocity $U = 0.25 \text{ m.s}^{-1}$. When the dipole is in an infinite domain and vanishing viscosity, it travels with a velocity $U$ along the direction $(\cos (\alpha), \sin(\alpha))$. In our case, the boundaries feature stress-free walls. The velocity perpendicular to the walls is zero, while the tangential velocity is indeterminate. When the dipole reaches the boundary, it splits into two vortices of opposite sign, which slide away from each other along the wall. In addition, the introduction of viscosity reduces the induced velocity of the dipole by homogenizing the intensity through the PSE scheme.
Because this problem does not have an explicit solution on a closed domain, we simulate the ground truth with the vortex method for fine particle discretization. The evolution of the ground truth is shown in Figure~\ref{fig:ref_trajectory}.

\begin{figure}[htbp]
    \centering
    \begin{subfigure}[b]{0.33\textwidth}
        \includegraphics[width=\textwidth]{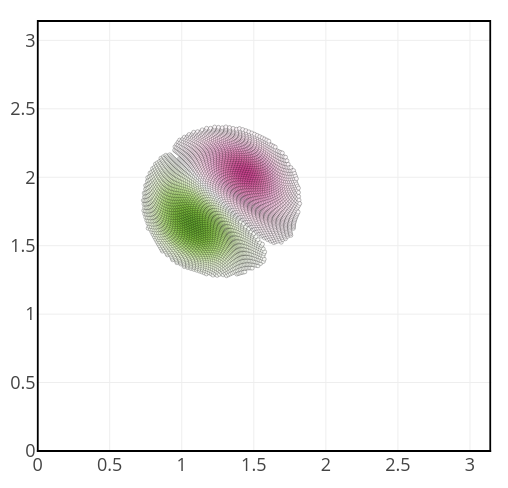}
    \end{subfigure}%
    \hfill
    \begin{subfigure}[b]{0.34\textwidth}
        \includegraphics[width=\textwidth]{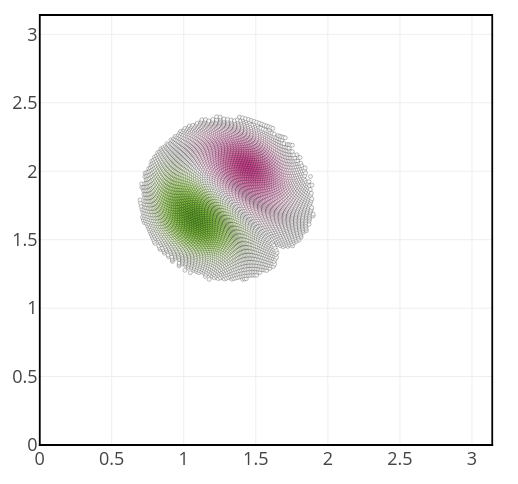}
    \end{subfigure}%
    \hfill
    \begin{subfigure}[b]{0.33\textwidth}
        \includegraphics[width=\textwidth]{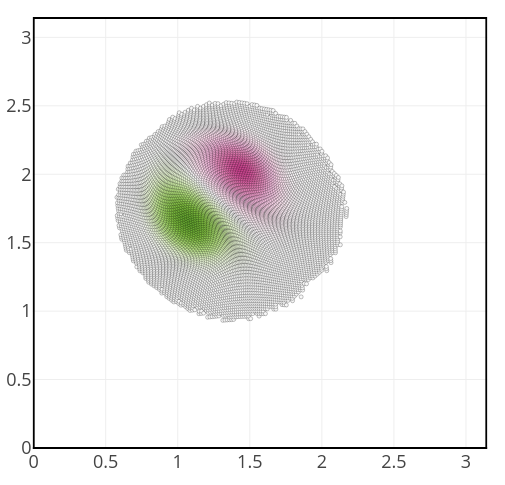}
    \end{subfigure}%
    \caption{Effect of the parameter $\varepsilon_\omega$ on the spatial extend of the particle discretization for one member. From left to right, results for $\varepsilon_\omega = 0.1, 0.01$, and $1.0e^{-6}$.}
    \label{fig:eps_effect}
\end{figure}

The vorticity field of the initial member is first discretized on a regular grid of particles with a characteristic length $d_p$. Each particle is assigned a circulation $\Gamma_p = \omega(\bm x_p) V_p$ where $V_p = d_p^2$.
Several numerical parameters in the simulation control the particle discretization and impact the solution accuracy. The first one is the particle size defined by $d_p$. Another significant parameter is the threshold value $\varepsilon_\omega$ associated with the remeshing processes during the forecast (to prevent high distortion of the particle distribution) and during the Remesh-EnKF filter. The threshold $\varepsilon_\omega$ is used to remove from the particle set all particles with intensity $|\Gamma_p| < \varepsilon_\omega V_p$. The impact of this parameter on the spatial extent of the particle distribution is illustrated in Figure~\ref{fig:eps_effect}.

\subsubsection{Assimilation settings}
An ensemble of $N=32$ members is created by sampling the dipole parameters' distributions. We sample independently the radius $R$, the prescribed velocity $U$, the orientation $\alpha$, and the center $\bm x_{\text{mean}}$. Additionally, the field viscosity $\nu$ is also sampled. The viscosity is not calibrated in these experiments but sampled to introduce further model error. The distributions of initial dipole characteristics and viscosity are reported in Table~\ref{tab:ens_dipole}.
\begin{table}[htbp]
    \centering
    \caption{Ensemble generation variables}
    \begin{tabular}[t]{|l|l|}
        \hline
        Variables     & Distributions                                                                    \\
        \hline
        radius        & $R \sim \cN(0.5, 0.05^2)$                                                        \\
        orientation   & $\alpha \sim \cU\left(\frac\pi2, \pi\right) (\text{rad.}) $                      \\
        dipole center & $\bm x_{\text{mean}} \sim \cN\left([\frac\pi2, \frac\pi2] ,0.1^2 \bm I_2\right)$ \\
        velocity      & $U \sim \cU(0.25, 0.5) $                                                         \\
        viscosity     & $v \sim \cN(0.0015, 0.0005^2)$                                                   \\
        \hline
    \end{tabular}
    \label{tab:ens_dipole}
\end{table}

\begin{figure}[htbp]
    \centering
    \makebox[1.05\textwidth][c]{%
        \includegraphics[width=1.05\textwidth]{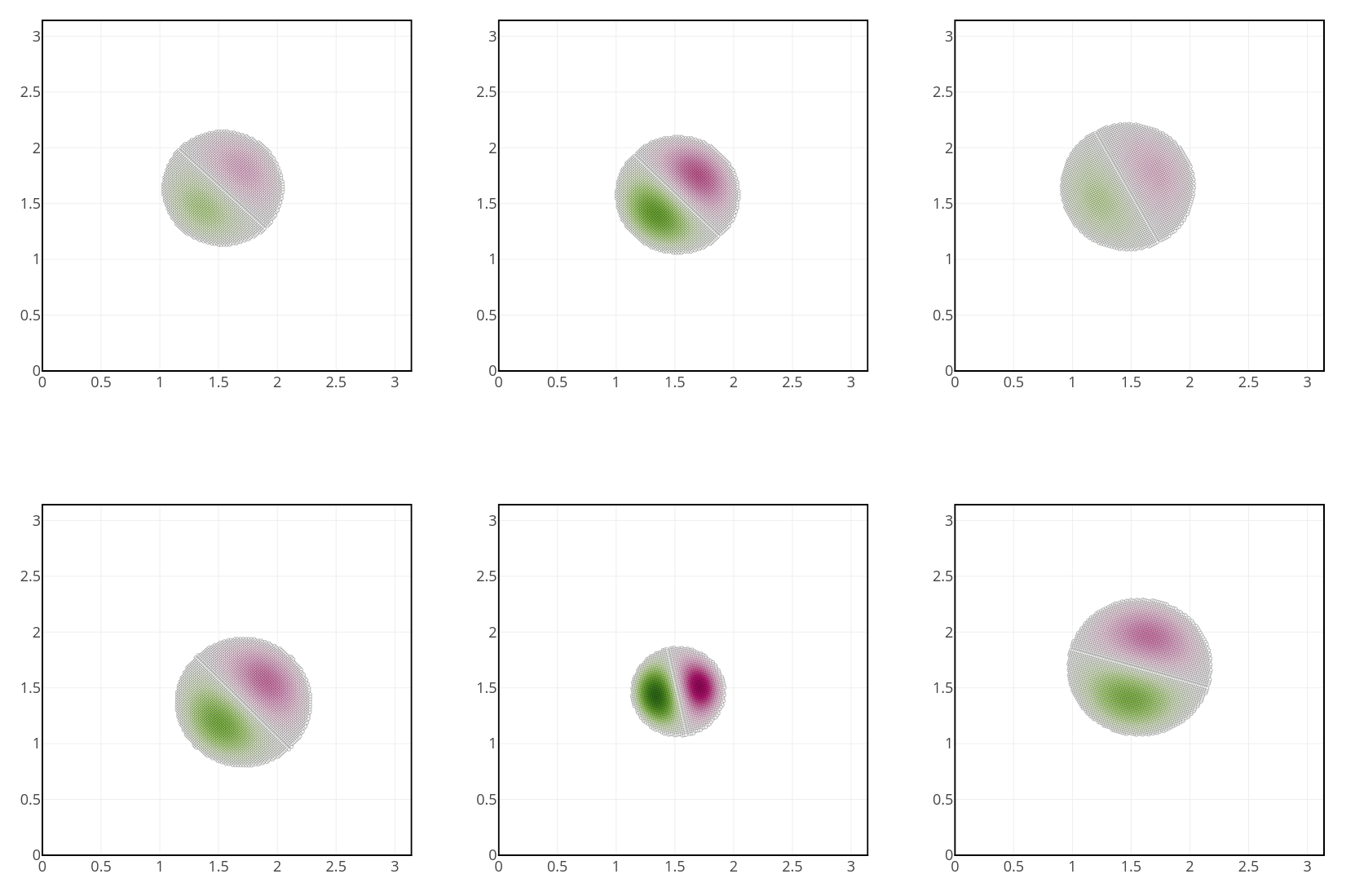}}
    \caption{Six samples from the initial ensemble.}
    \label{fig:sample_ens}
\end{figure}

Figure~\ref{fig:sample_ens} shows a sample set of the initial conditions.

\begin{figure}[htbp]
    \centering
    \includegraphics[width=0.9\textwidth]{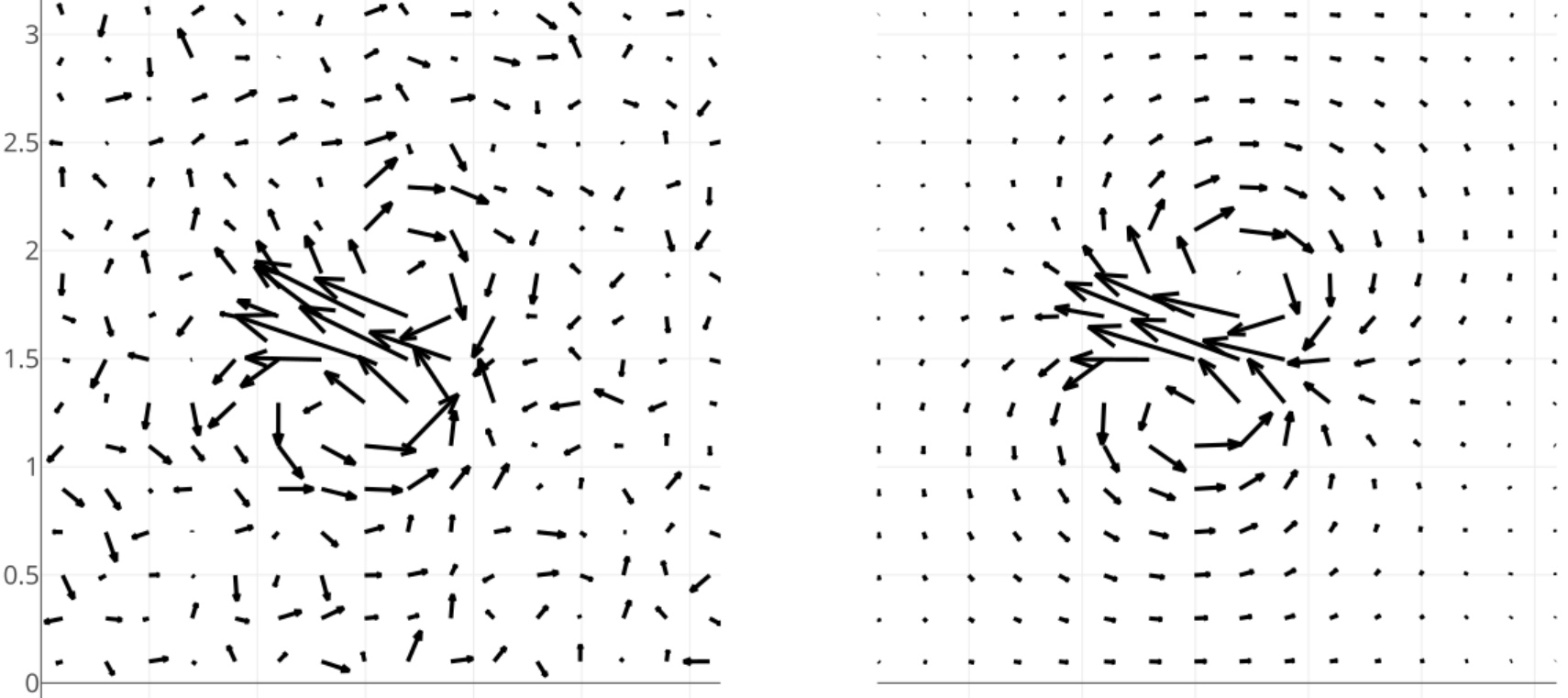}
    \caption{Observed and reference velocity fields. The error on each component is a sample from a centered normal distribution with the nominal value $\sigma_{\text{obs}} = 0.05$.}
    \label{fig:velocity}
\end{figure}

The assimilation frequency is defined by the assimilation step $\Delta t_a = 1$. Observations are collected on a regular grid of observation points of size $N_{\text{obs}} = 12^2$, measuring both components of the velocity field. The observations noise follows a normal distribution $\mathcal N(0, \sigma_{\text{obs}}^2 \bm{I})$, indicating an ensemble of independent measurements, each characterized by a standard distribution of $\sigma_{\text{obs}}$. An example of observed velocity with and without noise is shown in Figure~\ref{fig:velocity}. Table~\ref{tab:simu_2d} reports the assimilation and numerical parameters.
\begin{table}[htbp]
    \centering
    \caption{Nominal assimilation and simulation parameters}
    \begin{tabular}[t]{|l|l|}
        \hline
        Parameters                      & Values                           \\
        \hline
        time step                       & $dt = 0.005$                     \\
        final time                      & $t_f =10$                        \\
        std. observation                & $\sigma_{obs} = 0.05$            \\
        vorticity threshold             & $\varepsilon_{\omega} = 10^{-4}$ \\
        particle characteristic size    & $d_p = \frac{\pi}{256}$          \\
        smoothing length                & $\varepsilon = 2.0 d_p$          \\
        number of assimilation          & $N_{\text{assim}} = 10$          \\
        ensemble size                   & $N = 32$                         \\
        number of observation           & $N_{\text{obs}} = 12^2 $         \\
        grid discretization             & $N_{\text{grid}} = 128^2$        \\
        number of remeshing by forecast & $N_{\text{remesh}} = 2 $         \\
        \hline
    \end{tabular}
    \label{tab:simu_2d}
\end{table}

\subsubsection{Error definition}
We use an ensemble average \(L_2\)-error defined as
\begin{equation*}
    e_\omega = \frac1\nens \sum_{i = 1}^{\nens} \int_\Omega \left(\omega_i(\bm x) - \omega^{gt}(\bm x)\right)^2 \mathrm{d}\bm x
\end{equation*}and use the member errors to evaluate the dispersion of the error estimate. In practice, the integrals in the $L_2$ norms are computed by projecting the field onto a fine uniform grid, which is used for numerical quadrature.

\subsection{Results}
\subsubsection{Error evolution in time}
\begin{figure}[htbp]
    \centering
    \includegraphics*[width=0.9\textwidth]{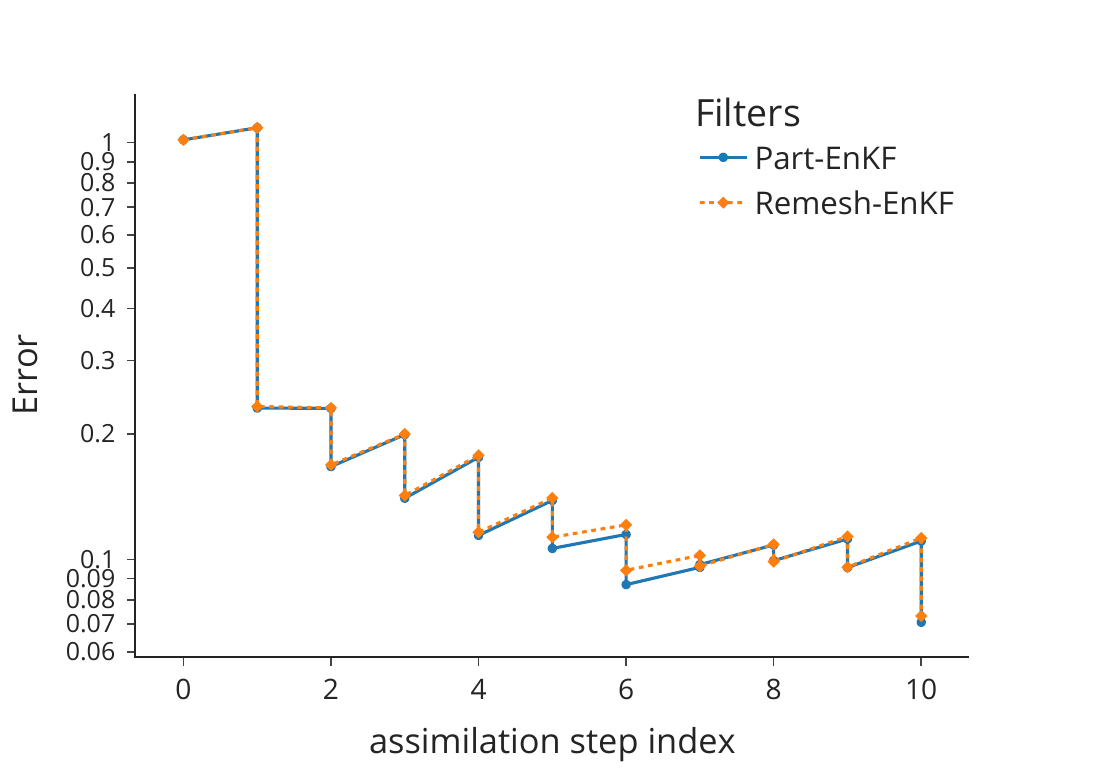}
    \caption{\(L_2\)-error curves of the field through assimilation steps.}
    \label{fig:assim_time}
\end{figure}

We start by analyzing the assimilation error over time. Figure~\ref{fig:assim_time} reports the error throughout the assimilation process for the two filters. Comparable error curves are reported for the two filters. At each assimilation step, the error decreases and slowly increases again until the next assimilation. Calibrating the viscosity would help reduce the growth of the error between 2 assimilation steps.

\subsubsection{Impact of assimilation parameters}
\begin{figure}[htbp]
    \centering
    \includegraphics[width=0.9\textwidth]{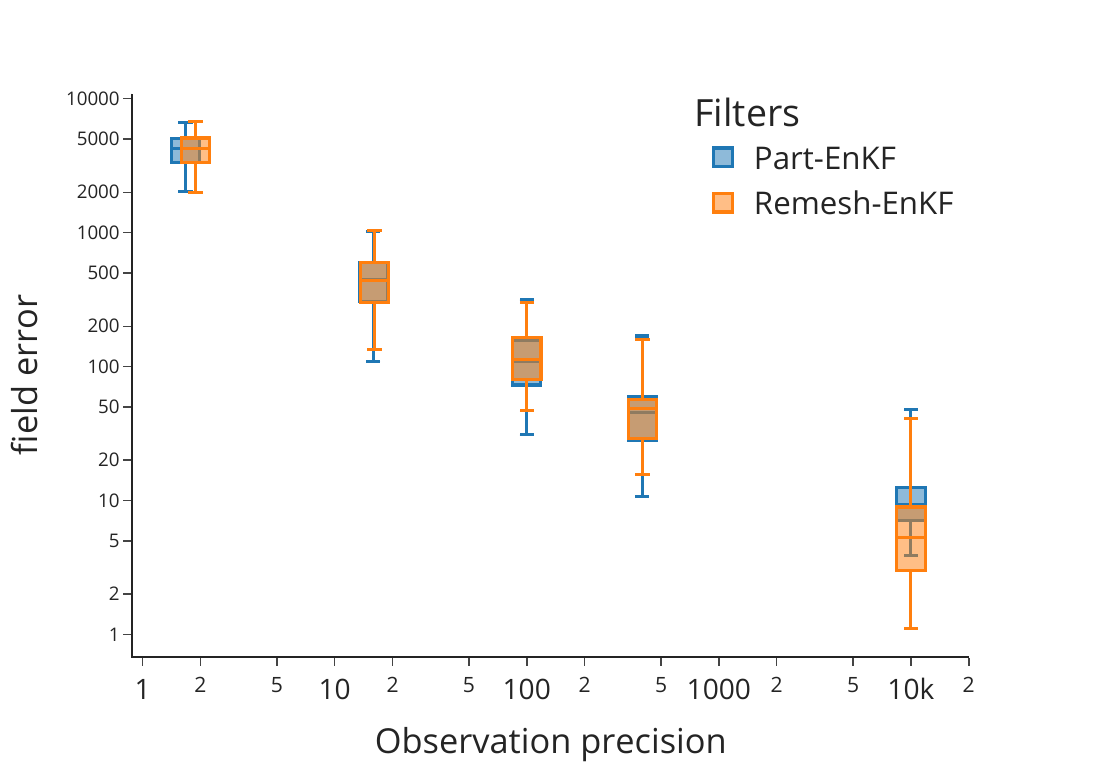}
    \caption{Box plots of the vorticity error w.r.t. observation prec. $1/\sigma_{\text{obs}}^2$.}~\label{fig:obs_precision_2}
\end{figure}

We assess the errors of the two filters for varying the assimilation parameters.
We observe the reduction of the assimilation error when the observation precision \(1/\sigma_{\text{obs}}^2\) increases, the number of observations \(N_{\text{obs}}\), and the number of assimilation step \(N_{\text{assim}}\).
Figure~\ref{fig:obs_precision_2}~shows that both filters exhibit a similar decrease in error bias and variance as observation precision improves, with a consistent convergence rate. The order of convergence is about 0.68 for Part-EnKF and 0.75 for Remesh-EnKF.

\begin{figure}[htbp]
    \centering
    \begin{subfigure}{0.9\textwidth}
        \centering
        \includegraphics[width=\textwidth]{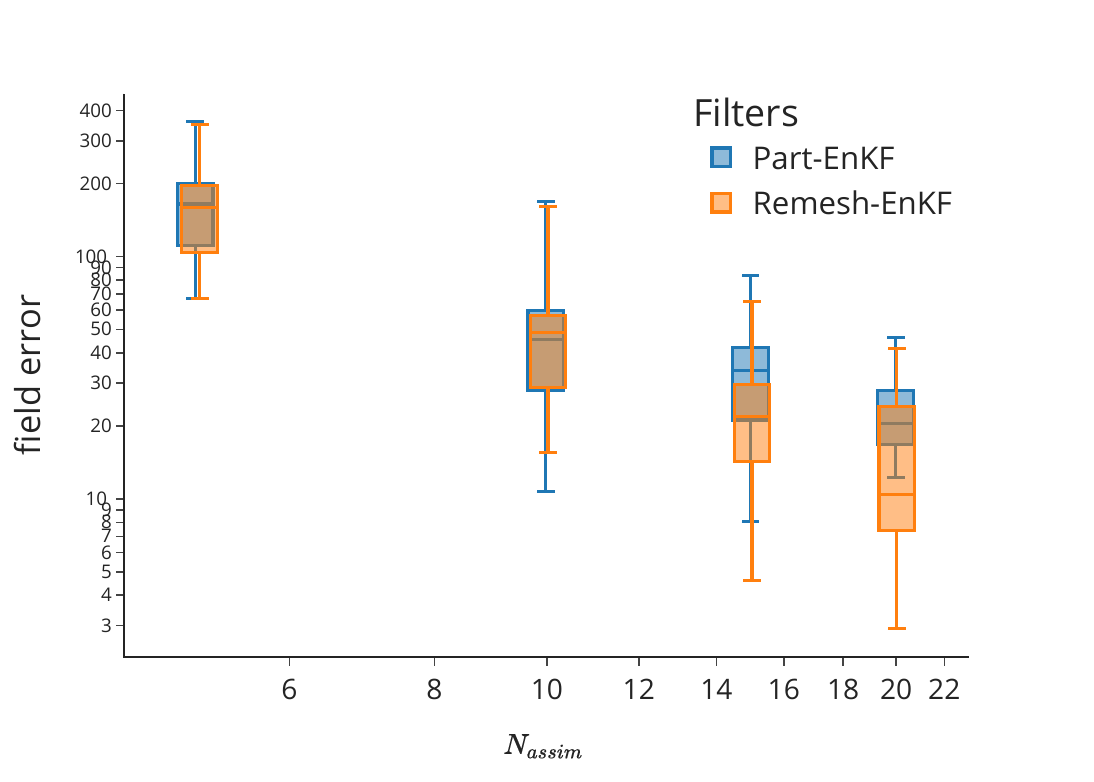}
    \end{subfigure}
    \caption{Box plots of the vorticity error w.r.t. $N_{\text{assim}}$.}~\label{fig:na}

\end{figure}

In Figure~\ref{fig:na}, the reduction of error is still prominent and shows a reduction of variance as the assimilation of the observation frequency increases. The error also decreases at a constant rate for both filters, with the Remesh-EnKF achieving a higher order of 1.8 compared to the Part-EnKF, which has an order of 1.4.

\begin{figure}[htbp]
    \centering
    \begin{subfigure}{0.9\textwidth}
        \centering
        \includegraphics[width=\textwidth]{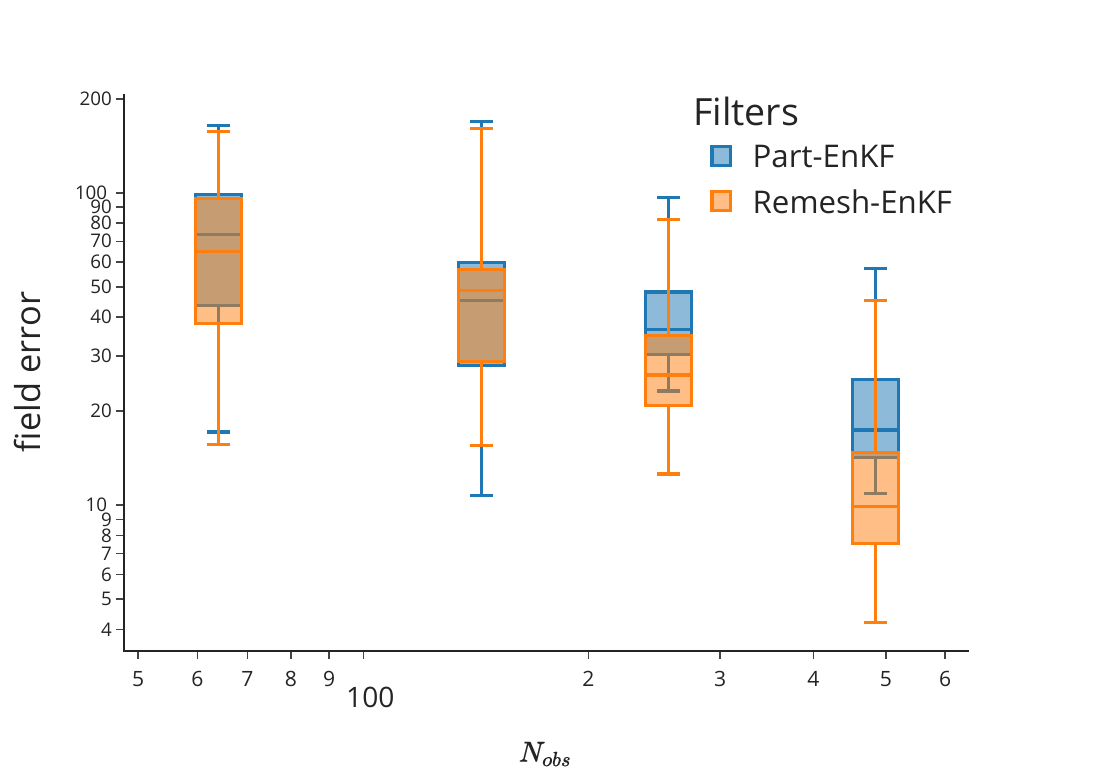}
    \end{subfigure}
    \caption{Box plots of the vorticity error w.r.t. $N_{\text{obs}}$.}~\label{fig:nobs}
\end{figure}

Finally, we analyze the error convergence with respect to the number of observations in Figure~\ref{fig:nobs}. Again, the bias and variances decrease when $N_{\text{obs}}$ increases. For a relatively small number of observations, both filters produce similar results, while the Remesh EnKF produces relatively better results as the number of observations becomes large. Moreover, the convergence rate seems to change around 200 observation points. Nevertheless, it illustrates adequate performances for both filters.

\subsubsection{Error with respect to simulation parameters}
\begin{figure}[htbp]
    \centering
    \begin{subfigure}[t]{0.32\textwidth}
        \centering
        \includegraphics[width=\textwidth]{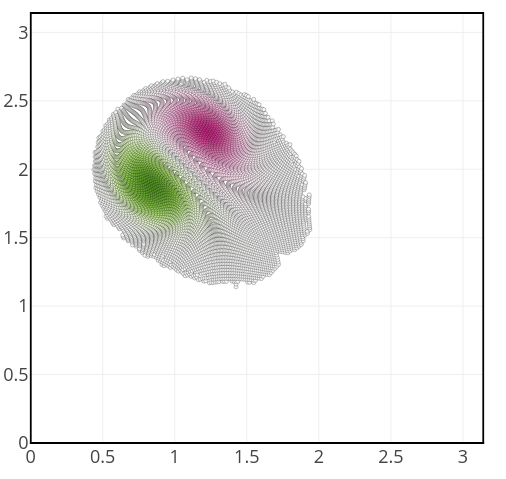}
        \caption{Forecast member.}
    \end{subfigure}
    \hfill
    \begin{subfigure}[t]{0.32\textwidth}
        \centering
        \includegraphics[width=\textwidth]{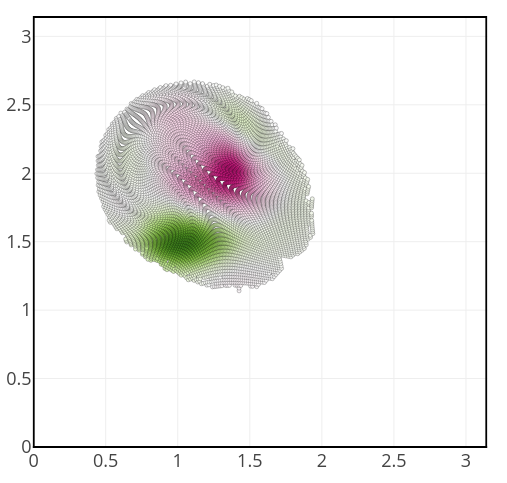}
        \caption{Part-EnKF analyzed solution.}
    \end{subfigure}
    \hfill
    \begin{subfigure}[t]{0.32\textwidth}
        \centering
        \includegraphics[width=\textwidth]{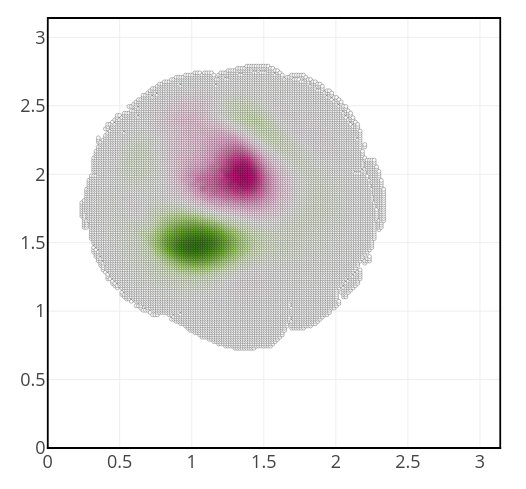}
        \caption{Remesh-EnKF analyzed solution.}
    \end{subfigure}
    \caption{Assimilation of one member with a forecast discretization unadapted to the analyses solution. The forecast discretization used by the Part-EnKF does not properly cover the support of the analyzed solution, introducing interpolation errors.}
    \label{fig:assim_member}
\end{figure}

To better understand the differences between the two filters, let us now turn to the dependence of the error on the particle discretization parameters. For the Part-EnKF, remember that each member has its own particle discretization, which is preserved by the assimilation. However, this scheme introduces different sources of error. First, irregularity in the particle distribution can cause approximation errors in the approximation of the analyzed solution. Even more critical, the member's particles may not be appropriately distributed over the whole support of the analyzed field they have to approximate. This effect can be seen in several members of the ensemble where the analyzed field is projected onto a nonconforming particle discretization. For instance, we detail the first assimilation step of one member for the different filters. The analyzed field is known over the entire spatial domain; we observed in Figure~\ref{fig:assim_member} that the Remesh-Filter is able to interpolate the entire solution. For Part-EnKF, the interpolation is limited to the support of the member's particles and is therefore incomplete. These effects are more noticeable when the assimilation frequency is low, resulting in higher errors, or when the thresholding is more aggressive. In addition, the approximations of Equation~\eqref{eq:first_part_approx} will introduce an error that increases for larger particle sizes $d_p$. In particular, unlike the Remesh EnKF, the Partial EnKF filter generally does not preserve any moment of the continuous analyzed field.

\begin{figure}[htbp]
    \centering
    \includegraphics[width=0.9\textwidth]{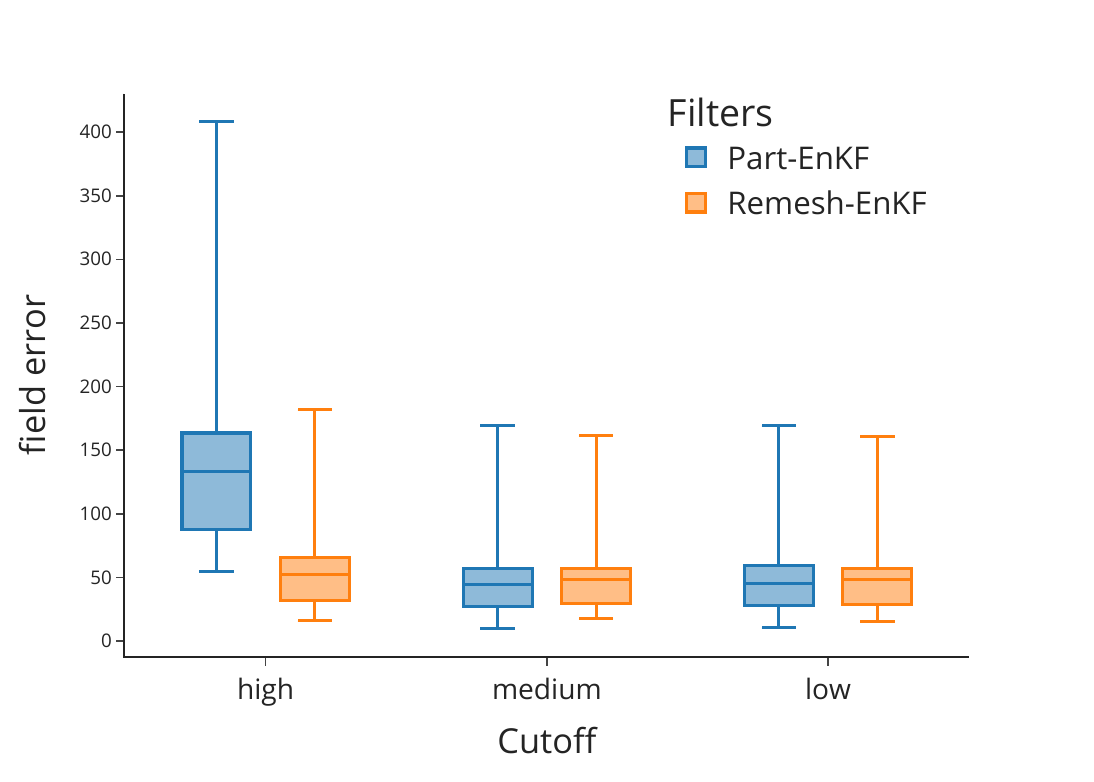}
    \caption{Vorticity error w.r.t. the threshold $\varepsilon_{\omega}$. The high, medium, and low cutoff correspond respectively to $\varepsilon_{\omega} = 0.1, 0.01$ and $1.e^{-6}$.}
    \label{fig:cuttoff}
\end{figure}

To investigate the effect of support in the particles' discretization, we varied the value of thresholding parameter $\varepsilon_{\omega}$. We have seen in Figure~\ref{fig:eps_effect} that this parameter affects the number of particles and, thus, the size of the support. In Figure~\ref{fig:cuttoff}, we observe that while Remesh-EnKF is not very sensitive to $\varepsilon_{\omega}$, Part-EnKF shows a much significant dependence of the error. 
\begin{figure}[htbp]
    \centering
    \includegraphics[width=0.9\textwidth]{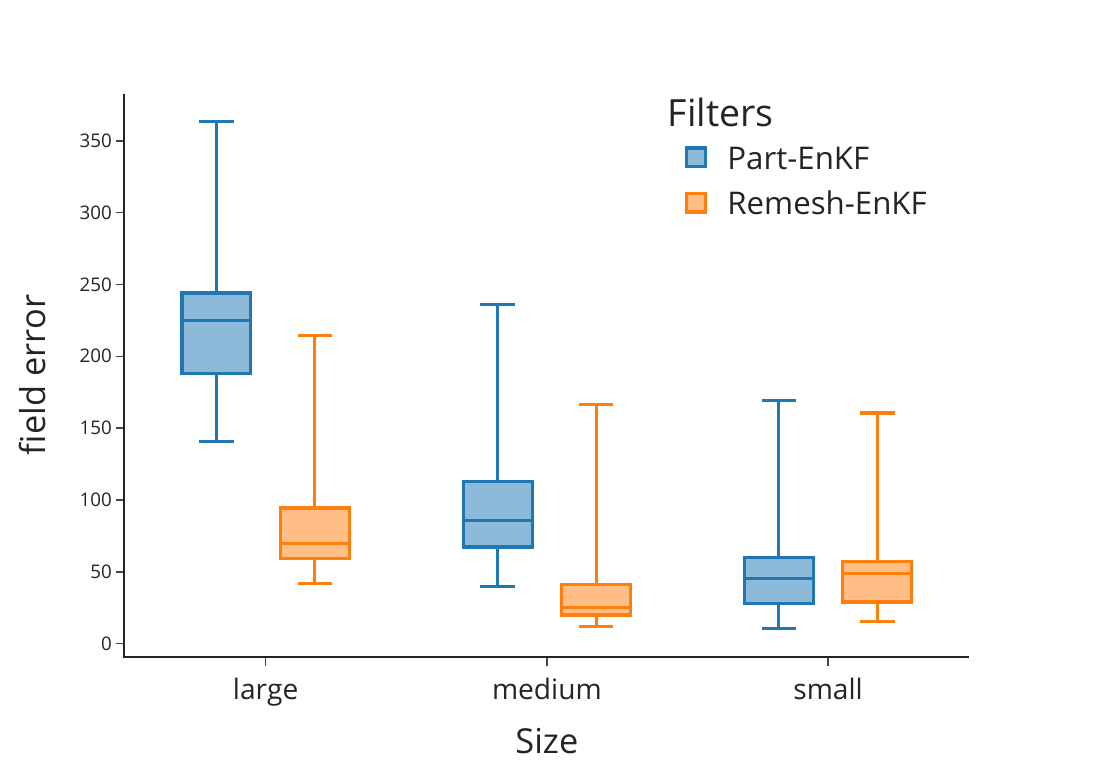}
    \caption{Vorticity error w.r.t. the particle size $d_p$. The large, medium, and low sizes correspond to $d_p = 0.0327, 0.0245$, and $0.0123$ respectively. The Part-EnKF error is strongly impacted by $d_p$ when Remesh-EnKF is less affected.}
    \label{fig:np}
\end{figure}

Concerning the particle size $d_p$, we evaluate its effect on the assimilation. In this section, we have used the first approximation defined in Equation~\eqref{eq:first_part_approx} to determine the new intensities for the Part-EnKF.This choice is motivated in order to accelerate the update process as it does not require system inversions or any iterations.
In Figure~\ref{fig:np}, we observe that the error for Part-EnKF increases with $d_p$ as expected. The behavior of the error confirms, in that case, the higher effect of the particle approximation for Part-EnKF. In contrast, the Remesh EnKF error is much less affected by much $d_p$. Using a regression operator in Part-EnKF will not solve the issue of inadequate coverage of the analyzed solution. This observation has motivated the development of the next Chapter.
\clearpage


\section{Conclusion}\label{sec:conclusion_chap4}
In this chapter, we introduced two methods for performing a sequential ensemble data assimilation approach with particle-based models. Specifically, we have developed two novel Ensemble Kalman Filter schemes dedicated to mesh-free Lagrangian simulations. These two methods rely either on updating the particle intensities by approximation or on remeshing the particle discretization on a common new grid of particles.

We proposed two update strategies, both relying on a correction matrix $\bm F$ that does not directly depend on the state discretization. The first strategy, Remesh-EnKF, is based on projecting all members onto a common set of particles for the discretization of the analyzed solution. The second strategy, Part-EnKF, evaluates the analyzed field at the members' particle locations to update their particle intensities without modifying the discretization.

The two methods have been first tested on a one-dimensional example and compared with a purely Eulerian method. The results demonstrate comparable performance provided that Part-EnKF uses sufficiently refined particle discretization so that the support of the particles is consistent with the analysis solution. However, in cases where the support deviates from the analysis field, members may have significant errors in Part-EnKF. Increasing the support for the solution is necessary.
We have also tested a two-dimensional case for a nonlinear Navier-Stokes equation. For various parameters analyzed, the two filters agree well.
As for the 1D case, the Part-EnKF yields a more significant error if the discretization support of the members does not overlap sufficiently. This effect is mitigated by extending the discretized support. In general, in situations where remeshing is not an option, it is necessary to improve the robustness of Part-EnKF. One possibility is to update not only the intensities but also to adjust the particle positions accordingly. Such adaptations could be derived from optimal transport schemes~\cite{bocquet_bridging_2023} or by aligning the background ensemble with observations through displacement correction methods~\cite{ravela_data_2007,rosenthal_displacement_2017}, as we propose in the next chapter.


\appendix
\bibliography{../../memoire/biblio}

\end{document}